\documentclass[12pt]{article}
\usepackage[german]{babel}
\usepackage{amssymb}
\usepackage{amsmath}
\usepackage{graphicx}
\usepackage{bbm}
\usepackage[right]{eurosym}
\usepackage[left=2cm,right=2cm,top=2cm,bottom=3cm]{geometry}
\def\vs{\vspace}
\def\noi{\noindent}

\def\IN{\mathbb N}
\def\IZ{\mathbb Z}
\def\IR{\mathbb R}

\def\IQ{\mathbb Q}
\def\IP{\mathbb P}

\def\an{\mathrm{an}}
\def\exp{\mathrm{exp}}

\def\ma{\mathcal}

\begin{document}
	\begin{center}
		{\bf \Large Lebesgue measure and integration theory on non-archimedean real closed fields with archimedean value group}
	\end{center}
	
	\centerline{Tobias Kaiser}
	
	\vspace{0.7cm}\noi \footnotesize {{\bf Abstract.} 
		Given a non-archimedean real closed field with archimedean value group which contains the reals, we establish for the category of semialgebraic sets and functions a full Lebesgue measure and integration theory such that the main results from the classical setting hold. The construction involves methods from model theory, o-minimal geometry and valuation theory. We set up the construction in such a way that it is determined by a section of the valuation. If the value group is isomorphic to the group of rational numbers the construction is uniquely determined up to isomorphism.
		The range of the measure and integration is obtained in a controlled and tame way from the real closed field we start with.
		The main example is given by the case of the field of Puiseux series where the range  is the polynomial ring in one variable over this field.

\section*{Introduction}
\label{Introduction}

\noi In the last two decades an abstract integration theory in algebraic and non-archimedean geometry, named motivic integration, has been set up to solve deep geometric problems.
It has been introduced by Kontsevich in 1995 and was further developed by Denef and Loeser [18], by Cluckers and Loeser [8, 9] and others to obtain, using geometry, valuation theory and model theory, a measure and integration theory on Henselian discretely valued fields (for example $p$-adic integrals). Next, Hrushovski and Kazhdan [29] have defined motivic integration for algebraically closed valued fields and Yimu Yin [42] has developed an analogue for o-minimal valued fields.
The ranges of motivic integration are abstract spaces as Grothendieck rings of definable sets.
We refer to Cluckers et al. [13, 14] for an overview on motivic integration, especially to the introduction of their Volume I.

\vs{0.1cm}
\noi We develop a measure and integration theory in Lebesgue's style for the setting of ordered  fields. Since every ordered field distinct from the field of reals is totally disconnected and not locally compact one has to restrict to a tame setting to have a chance to develop reasonable analysis.
There exists a nice differential calculus in semialgebraic geometry over arbitrary real closed fields (see Bochnak et al. [3, 2.9]) or, more general, in o-minimal structures over (necessarily real closed) fields (see Van den Dries [19, Chapter 7]).
For semialgebraic or o-minimal geometry over the reals the classical Lebesgue measure and Lebesgue integration is available and has been successfully used for geometric questions (see for example Yomdin and Comte [43]).
One can also establish easily a Lebesgue measure and integration theory for an archimedean real closed field by embedding the latter into the field of reals (the construction can be also performed via the classical limit process, see [31]).
But for general, i.e. non-archimedean, real closed fields such an integration theory does not exist so far.
There have been defined real valued measures (see Hrushovski et al. [30]) and integration of piecewise constant definable functions with respect to the Euler characteristic in a motivic style (see Br\"ocker [6] and Cluckers and Edmundo [7]), but these measures lack several analytic and geometric properties of the Lebesgue measure.
Preceeded by work of Berarducci and Otero [2] and of Fornasiero and Vasquez Rifo [26],
Ma\v{r}\'{i}kov\'{a} and Shiota [38, 39] have introduced via a limit process, mimicking the construction of the Lebesgue measure on the reals, a restricted measure for definable sets in o-minimal fields. It is defined for bounded sets and the range is just a semiring.
They obtain a partial transformation formula.
Recently, Costin, Ehrlich and Friedman [16] have developed an integration theory on the surreals in the univariate case. They also show that a reasonable integration theory in the non-archimedean setting requires in general a tame setting.

 \vs{0.75cm}
\hrule

\vs{0.4cm}
{\footnotesize{\itshape 2010 Mathematics Subject Classification:} 03C64, 03H05, 06F20, 12J25, 14P10, 28B15, 28E05, 32B20}
\newline
{\footnotesize{\itshape Keywords and phrases:} Non-archimedean real closed field, archimedean value group, semialgebraic sets and functions, Lebesgue measure and integration}
\newline
{\footnotesize{\itshape Acknowledgements:} The authors were supported in parts by DFG KA 3297/1-2.}

\newpage
\noi Given a non-archimedean real closed field with archimedean value group which contains the reals, we are able to construct a full Lebesgue measure and Lebesgue integration theory for the category of semialgebraic sets and functions with good control over the range such that the main properties of the real Lebesgue measure and integration hold:\\
The Lebesgue measure for semialgebraic sets on such a field is
finitely additive, monotone, translation invariant and
reflects elementary geometry.
The latter means, for instance, that the measure of an interval is, as it should be, the length of the interval.\\
The Lebesgue integral in this setting is linear,
the transformation formula,
Lebesgue's theorem on dominated convergence and the fundamental theorem of calculus hold with the necessary adjustments.
(For example, for Lebesgue's theorem we have to work, as usual in semialgebraic geometry, with one-parameter families instead of sequences.)
Moreover, a version of Fubini's theorem can be established.\\
Our construction relies on results and methods from model theory, o-minimal geometry and valuation theory.

\vs{0.1cm}
\noi First we want to mention that such a Lebesgue measure and integration theory cannot in general be performed inside a given real closed field, the integrals take necessarily values outside the given field. This can be seen by looking at the function $x\mapsto \int_1^x dt/t$ for $x>0$. A very basic version of the transformation formula gives that this function defines a logarithm. But on certain real closed fields there are no reasonable logarithms with values in the field (see Kuhlmann et al. [35]). Moreover, assuming very basic facts for the measure as the above property that the measure of an interval is its length, it cannot be $\sigma$-additive if the range is an ordered ring extension (consider the cut of the given field obtained by the natural numbers and the infinitely large elements). So a construction via a limit process is not possible.\\
We show that, given a non-archimedean real closed field $R$ with archimedean value group which contains the reals, it is enough to work in the immediate maximal extension (see for example Kaplansky [34]) with respect to the standard real valuation and then to add a logarithm to obtain the range for a reasonable measure and integration theory.

\vs{0.1cm}
\noi The basic idea of the construction of the measure (and similarly of the integral) is the following:\\
Let $A$ be a semialgebraic subset of some $R^n$. Then $A$ can be expressed by a formula in the language of ordered rings augmented by symbols for the elements of $R$ (by Tarski even by finitely many equalities and inequalities of polynomials over $R$). Replacing the tuple $a$ of elements from $R$ involved in this formula by variables, we get a parameterized semialgebraic family of subsets of $\IR^n$.
For this family we can apply the results of Comte, Lion and Rolin [15, 37] (see also the related work of Cluckers and D. Miller [10, 11, 12]) and obtain that the function $F$ given by computing the Lebesgue measure of the members of the family is given by a polynomial in globally subanalytic functions depending on the parameter and their logarithms (see Van den Dries and Miller [24] for the notion of globally subanalytic sets and functions).
Now we want to plug in the above tuple $a$ in this function and to declare this as the Lebesgue measure of the given set $A$.
For that purpose we need a suitable extension of $R$ where this can be done in a reasonable way. One could think of some ultrapower of $\IR$ containing $R$.
Here we would enter non-standard measure and integration theory (see for example Robinson [41] and Cutland [17]). But by this abstract choice we would loose control about the values obtained by measuring and integrating.
Instead, we use the work of Van den Dries, Macintyre and Marker [21, 22, 23] and choose the following embedding, carried out in two steps. Let $\Gamma$ be the value group of $R$ with respect to the standard real valuation.
Then we embed $R$ into the power series field $\IR((t^\Gamma))$. This is a model of the theory $T_\an$ of the o-minimal structure $\IR_{\an}$ (see [24]; the sets definable in $\IR_\an$ are exactly the globally subanalytic sets), so that we can evaluate the globally subanalytic functions at the tuple $a$.
Since $\Gamma$ is archimedean there is an order preserving embedding of $\Gamma$ in the additive group of the reals.
Hence we can view $\IR((t^\Gamma))$ as a subfield of $\IR((t^\IR))$. The latter is finally embedded in the field $\IR((t))^{\mathrm{LE}}$ of logarithmic-exponential series (in short LE-series). From the description of $F$ and the definition of the logarithm on power series we see that $F(a)$ lies in the ring $\IR((t^\Gamma))\big[\log(t^\Gamma)]$. With $\Gamma$ viewed as a subgroup of $\IR$, the properties of the logarithm on the field of LE-series give that
the latter ring equals $\IR((t^\Gamma))\big[X\big]$ where $X:=\log(t^{-1})$ is independent over $\IR((t^\Gamma))$. We call this polynomial algebra over $\IR((t^\Gamma))$ the
Lebesgue algebra of $R$.\\
Every such choice of an embedding of the given real closed field into the field $\IR((t^{\mathrm{LE}}))$ gives a Lebesgue measure and integral (with all the desirable properties). To compare the construction obtained by different choices of embeddings, we introduce a natural notion of equivalence. We call the Lebesgue measure, respectively integral, with respect to different embeddings isomorphic if they are the same up to an isomorphism of the range, i.e. of the Lebesgue algebra.
We prove that we obtain isomorphy up to the choice of a section for $R$; i.e. a group homomorphism from the value group of $R$ into the multiplicative group of its positive elements that is compatible with the valuation. Hence we can describe the Lebesgue measure and integral in internal terms of the given real closed field.
Moreover, if the value group is given by the rationals then the construction is, up to isomorphisms, uniquely determined. In the case that the given field is the field $\IP$ of real Puiseux series or convergent real Puiseux series, we obtain a nice restriction of the range of the measure and the integral; it is given by the polynomial ring $\IP[X]$ in one variable over $\IP$.

\vs{0.1cm}
\noi We have explained how we construct the measure (and similarly the integral). Now we explain why the above mentioned results from classical Lebesgue measure theory and integration theory hold in our setting, independently from the chosen section. The above function $F$ is definable in the o-minimal structure $\IR_{\an,\exp}$ (see [21, 24] for this important structure).
We have embedded $R$ into the non-standard model $\IR((t^{\mathrm{LE}}))$ of the theory $T_{\an,\exp}$ of $\IR_{\an,\exp}$ to plug in the tuple $a$.
We formulate the properties of classical Lebesgue measure and integration theory as statements in the natural language of $\IR_{\an,\exp}$ and can then transfer these results to this non-standard model. But we want to obtain the results for raw data from $R$. This reduction requires some work using the set-up of the construction and o-minimality.

\vs{0.1cm}
\noi The condition that the given real closed field $R$ contains the reals was introduced to keep the technical details at a reasonable level. Note that the real numbers can be adjoined in a unique way to an arbitrary real closed field and that the value group stays the same (see Prie\ss-Crampe [40, III \S 1]).

\vs{0.1cm}
\noi The paper is organized as follows.
After introducing the notations and terminology used throughout the work, we present and develop in Section 1 the setting and background for our constructions and results.
We cover archimedean ordered abelian groups, the standard real valuation on real closed fields, the theory of power series fields over the reals including the construction of the partial logarithm, the o-minimal structures $\IR_\an$ and $\IR_{\an,\exp}$, and, finally, the results on integration of parameterized definable functions.
Section 2 is devoted to the construction of the measure and the integral, and their elementary properties are shown.
The construction is performed with respect to a certain tuple of raw data called Lebesgue datum. We define in Section 3 a natural notion of equivalence between the results obtained from different Lebesgue data. We show that our construction does only depend on the choice of a section for the given real closed field. In the case that the value group is the group of rationals we obtain actually that the constructions are isomorphic. We present the main example given by the geometrically significant field of Puiseux series. Moreover, we consider extensions of real closed fields and the behaviour of the constructed measure with respect to the standard part map.\\
For the rest of the paper we develop our theory for the field $\IP$ of real Puiseux series to keep the notations at a reasonable level. The functions obtained by integrating are given by so-called constructible functions (compare with Cluckers and D. Miller [10, 11, 12]). We introduce in Section 4 these functions on $\IP$ which takes values in the polynomial ring $\IP[X]$ and develop the necessary analysis. This requires some work since we deal with a kind of hybrid of the $T_\an$-model $\IP$ and the $T_{\an,\exp}$-model $\IR((t^{\mathrm{LE}}))$.
We use these results in Section 5 to present the main theorems of integration, namely the transformation formula, Lebesgue's theorem on dominated convergence, the fundamental theorem of calculus and Fubini's theorem.
The final Section 6 is devoted to an application.
We exploit the smoothing properties of convolution to show an approximation result for unary functions definable in the $T_\an$-model $\IP$.

\vs{0.1cm}
\noi In upcoming papers we deal with the case of arbitrary non-archimedean real closed fields including the surreals, develop integration on semialgebraic manifolds including Stokes' theorem and apply the measure and integration theory on the field of Puiseux series to geometric questions on semialgebraic sets over the reals.

\section*{Notations}

\noi Throughout the paper we assume basic knowledge of the theory of real closed fields and semialgebraic geometry (see [3]), of o-minimal structures (see Van den Dries [19]), of model theory (see Hodges [28]) and of
measure and integration theory (see Bauer [1] and Bourbaki [4, 5]).

\vs{0.1cm}
\noi By $\IN=\big\{1,2,3,\ldots\big\}$ we denote the set of natural numbers and by $\IN_0=\big\{0,1,2,\ldots\big\}$ the set of natural numbers with $0$.

\noi Let $R$ be a real closed field. We set $R^*:=\{x\in R\mid x\neq 0\}, R_{>0}:=\{x\in R\mid x>0\}$ and $R_{\geq 0}:=\{x\in R\mid x\geq 0\}$. For $a,b\in R$ with $a<b$ let
$]a,b[_R=]a,b[:=\{x\in R\mid a<x<b\}$ be the open and $[a,b]_R=[a,b]:=\{x\in R\mid a\leq x\leq b\}$ be the closed interval with endpoints $a$ and $b$.
By $|x|$ we denote the euclidean norm of $x\in R^n$.

\noi Given a subset $A$ of $R^n$ we denote by $\mathbbm{1}_A$ the characteristic function of $A$.
For a function $f:R^n\to R$ we set $f_+:=\max(f,0)$ and $f_-:=\max(-f,0)$. By $\mathrm{graph}(f)$ we denote the graph of $f$. If $f$ is non-negative we call $\mathrm{subgraph}(f):=\big\{(x,s)\in R^{n+1}\mid 0\leq s\leq f(x)\big\}$ the subgraph of $f$.
For $A\subset R^{q+n}$ and $t\in R^q$ we set $A_t:=\{x\in R^n\mid (t,x)\in A\}$.
For $f:R^{q+n}\to R$ and $t\in R^q$ we define $f_t:R^n\to R, f_t(x)=f(t,x)$.
By $D_\varphi$ we denote the Jacobian of a partially differentiable function $\varphi:U\to R^n$ where $U\subset R^m$ is open (in the euclidean topology).

\noi Let $R[X]$ be the polynomial ring over $R$ in one variable. We equip the polynomial ring with the standard degree $\mathrm{deg}$ and set
$R[X]_{\leq n}:=\{f\in R[X]\mid \mathrm{deg}(f)\leq n\big\}$ for $n\in\IN_0$.

\noi Let $\ma{L}_{\mathrm{or}}=\big\{+,-,\cdot,<,0,1\big\}$ be the language of ordered rings. By Tarski, the theory of real closed fields has quantifier elimination in this language.
Given a formula $\varphi$ in the language, the notation $\varphi(x_1,\ldots,x_n)$ indicates that the free variables of $\varphi$ are among the variables $x_1,\ldots,x_n$.

\noi Given an extension $R\subset S$ of real closed fields and a semialgebraic subset $A$ of $R^n$ or a semialgebraic function $f:R^n\to R$ we denote by $A_S$ and $f_S$ the canonical lifting of $A$ and $f$ to $S$, respectively (see [3, Chapter 5]).

\noi The above notations are used analogously in other situations if applicable.

\noi Finally, by $\infty$ we denote an element that is bigger than every element of a given ordered set.

\section{Preparations}

\noi In this section we present and develop the necessary background for the constructions and results of the paper.

\subsection{Ordered abelian groups}

\noi An {\bf ordered abelian group} is an additively written abelian group with total ordering $\leq$ which is compatible with the addition.
We refer to Prie\ss-Crampe [40, Kapitel I] and Fuchs [27, Kapitel IV] for more information and the proofs of the results below.

\vs{0.2cm}
\noi Let $\Gamma$ be an ordered abelian group.
The group $\Gamma$ is {\bf divisible} if for every $\gamma\in\Gamma$ and $n\in\IN$ there is some $\delta\in\Gamma$ with $n\delta=\gamma$.
Divisibility holds if and only if $\Gamma$ is a $\IQ$-vector space.

\vs{0.2cm}
\noi The absolute value of $\gamma\in\Gamma$ is defined by $|\gamma|:=\max\{\gamma,-\gamma\}$.
The group $\Gamma$ is called {\bf archimedean} if for every $\gamma,\delta\in\Gamma\setminus\{0\}$ there is some $n\in\IN$ with $|\gamma|\leq n|\delta|$.

\vs{0.5cm}
{\bf 1.1 Fact} (H\"older)

\vs{0.1cm}
{\it Assume that $\Gamma$ is archimedean. Then there is an embedding $\Gamma\hookrightarrow (\IR,+)$ of ordered groups.}

\vs{0.5cm}
So the archimedean ordered abelian groups are, up to isomorphisms of ordered groups, precisely the subgroups of the additive group of $\IR$.\\
Moreover, given $\gamma\in \Gamma_{>0}$ there is a uniquely determined embedding $\varphi:\Gamma\hookrightarrow (\IR,+)$ of ordered groups with $\varphi(\gamma)=1$.

\subsection{Standard real valuation}

\noi Let $R$ be real closed field.
It is called {\bf archimedean} if the ordered group $(R,+)$ is archimedean.
We refer to [27, Kapitel VII] for the following.

\vs{0.5cm}
{\bf 1.2 Fact} (H\"older)

\vs{0.1cm}
{\it Assume that $R$ is archimedean. Then there is a unique field embedding $R\hookrightarrow \IR$.}

\vs{0.5cm}
\noi Let $R$ be an arbitrary real closed field.
The set
\[\ma{O}_{R}:=\big\{f\in R\mid -n\leq f\leq n\mbox{ for some }n\in \IN\big\}\]
of {\bf bounded} elements of $R$ is a valuation ring of $R$ with
maximal ideal
\[\mathfrak{m}_{R}:=\big\{f\in R\mid -1/n<f<1/n\mbox{ for all }n\in \IN\big\}\]
consisting of the {\bf infinitesimal} elements of $R$. Note that $\ma{O}_{R}$ is a convex subring of $R$.
Let $v_{R}:R^*\to R^*/\ma{O}_{R}^*$ be the corresponding valuation with value group $\Gamma_{R}:= R^*/\ma{O}_{R}^*$. It makes $R$ an ordered valued field, meaning that $0<f\leq g$ implies $v_{R}(f)\geq v_{R}(g)$. Note that the value group $\Gamma_{R}$ is divisible and that $\Gamma_R=\{0\}$ if and only if $R$ is archimedean.
We denote the residue field $\ma{O}_{R}/\mathfrak{m}_{R}$ by $\kappa_{R}$.
The residue field is an archimedean real closed field. If $R$ contains the reals then $\kappa_R=\IR$.

\vs{0.2cm}
\noi A {\bf section} for $R$ is a homomorphism from the the value group $\Gamma_R$ to the multiplicative group $R_{>0}$ such that $v_R(s(\gamma))=\gamma$ for all $\gamma\in \Gamma_R$. Since $\Gamma_R$ is divisible there is always a section for $R$.

\subsection{Power series fields over the reals}

\noi We refer to [40, Kap. II \S 5] and Van den Dries et al. [21, Section 1.2] for the following.

\vs{0.2cm}
\noi Let $\Gamma=(\Gamma,+)$ be an additively written ordered abelian group.
We consider the {\bf power series field} $\ma{R}:=\IR((t^\Gamma))$.
The elements of $\ma{R}$ are the formal power series $f=\sum_{\gamma\in \Gamma}a_\gamma t^\gamma$ with exponents $\gamma\in \Gamma$ and coefficients $a_\gamma\in \IR$ such that the support of $f$, $\mathrm{supp}(f):=\{\gamma\in\Gamma\mid a_\gamma\neq 0\}$, is a well-ordered subset of $\Gamma$.

\vs{0.2cm}
\noi {\bf Field structure:}
The addition given by
\[\big(\sum_{\gamma\in\Gamma}a_\gamma t^\gamma\big)+\big(\sum_{\gamma\in\Gamma}b_\gamma t^\gamma\big)=\sum_{\gamma\in\Gamma}(a_\gamma+b_\gamma)t^\gamma\]
and the multiplication given by
\[\big(\sum_{\gamma\in\Gamma}a_\gamma t^\gamma\big)\cdot\big(\sum_{\gamma\in\Gamma}b_\gamma t^\gamma\big)=\sum_{\gamma\in\Gamma}(\sum_{\alpha+\beta=\gamma}a_\alpha b_\beta)t^\gamma\]
establish a field structure on
$\ma{R}$ (note that for every $\gamma\in \Gamma$ the sum $\sum_{\alpha+\beta=\gamma}a_\alpha b_\beta$ is finite since the supports are well-ordered). In particular, $t^\gamma\cdot t^\delta=t^{\gamma+\delta}$ for all $\gamma,\delta\in \Gamma$. We identify
$\IR$ with a subfield of $\ma{R}$ by the field embedding $\IR\hookrightarrow \ma{R}, a\mapsto at^0$.

\vs{0.2cm}
\noi {\bf Ordering:}
By setting $\sum_{\gamma\in\Gamma}a_\gamma t^\gamma<\sum_{\gamma\in\Gamma}b_\gamma t^\gamma$ if $a_\delta<b_\delta$ where $\delta=\min\big\{\gamma\in\Gamma\mid a_\gamma\neq b_\gamma\big\}$, the field $\ma{R}$ becomes an ordered field.

\vs{0.2cm}
\noi The ordered field $\ma{R}$ is real closed if and only if $\Gamma$ is divisible. From now on we assume this.

\vs{0.2cm}
\noi {\bf Valuation:}
The valuation $v_{\ma{R}}$ is given by $\mathrm{ord}:\ma{R}^*\to \Gamma, f\mapsto \min \big(\mathrm{supp}(f)\big)$. We have $\Gamma_{\ma{R}}=\Gamma$,
\[\ma{O}_{\ma{R}}=\IR((t^{\Gamma_{\geq 0}})):=\big\{f\in \ma{R}\mid \mathrm{supp}(f)\subset \Gamma_{\geq 0}\big\}\]
and
\[\mathfrak{m}_{\ma{R}}=\IR((t^{\Gamma_{>0}})):=\big\{f\in \ma{R}\mid \mathrm{supp}(f)\subset \Gamma_{> 0}\big\}.\]
Note that $\kappa_{\ma{R}}=\IR$.

\vs{0.5cm}
{\bf 1.3 Fact}

\vs{0.1cm}
	{\it Let $R$ be a real closed field that contains the reals and let $\ma{R}:=\IR((t^{\Gamma_R}))$. Let $s:\Gamma_R\to R_{>0}$ be a section for $R$. Then there is a field embedding $\sigma:R\hookrightarrow \ma{R}$ such that $\sigma(s(\gamma))=t^\gamma$ for all $\gamma\in\Gamma_R$. Such an embedding is valuation and order preserving.}

\vs{0.5cm}
\noi We say that $\sigma$ is an embedding with respect to the section $s$. Note that it is in general not uniquely determined.
Note also that such an embedding is the same as having a valuation preserving embedding with $t^\Gamma$ in its image.

\vs{0.2cm}
\noi The power series field $\ma{R}=\IR((t^{\Gamma}))$ carries
a {\bf partial logarithm}
\[\log_\ma{R}=\log: \big(\IR_{>0}+\mathfrak{m}_\ma{R},\cdot\big)\stackrel{\cong}{\longrightarrow} \big(\ma{O}_\ma{R},+\big)\]
extending the logarithm on the reals
(compare with [21, Section 1.2] and S. Kuhlmann [36])
which is defined as follows:

\noi Let $f\in \IR_{>0}+\mathfrak{m}_\ma{R}$. Then there are unique $a\in\IR_{>0}$ and $h\in\mathfrak{m}_\ma{R}$ such that
$f=a(1+h)$. Then $\log(f)=\log(a)+L(h)$ where
\[L(x)=\sum_{j=1}^\infty \frac{(-1)^{j+1}}{j}\, x^j\]
is the {\bf logarithmic series}.

\noi The partial logarithm gives an order isomorphism between the multiplicative group of positive units of the ordered valuation ring $\ma{O}_\ma{R}$
and the additive group of the latter. Its inverse is given by the partial exponential function.

\subsection{The o-minimal structures $\IR_\an$ and $\IR_{\an,\exp}$}

\noi For $n\in\IN_0$ let $\IR\{x_1,\ldots,x_n\}$ denote the ring of convergent real power series in $n$ variables
and
\[\IR\langle x_1,\ldots,x_n\rangle:=\Big\{f\in \IR\{x_1,\ldots,x_n\}\,\Big\vert\, f\mbox{ converges on a neighbourhood of }[-1,1]^n\Big\}.\]
Note that we obtain in the case $n=0$ the real field $\IR$. Given $f\in \IR\langle x_1,\ldots,x_n\rangle$, the function
\[\widetilde{f}:\IR^n\to\IR, x\mapsto \left\{\begin{array}{lll}
f(x),&&x\in [-1,1]^n,\\
&\mbox{if}&\\
0,&&x\notin [-1,1]^n,\\
\end{array}\right.\]
is called a {\bf restricted analytic function}. (Note that $\IR^0=\{0\}$.)

\noi The language $\ma{L}_\an$ is obtained by augmenting the language  $\ma{L}_{\mathrm{or}}=\big\{+,-,\cdot,<,0,1\big\}$ of ordered rings by a function symbol for every restricted analytic function.
Let $\IR_\an$ be the natural $\ma{L}_\an$-structure on the real field. Sets resp. functions definable in $\IR_\an$ are the globally subanalytic sets resp. functions. These are the sets resp. functions that are subanalytic in the ambient projective space (see Van den Dries and Miller [24, p. 505]).
The $\ma{L}_\an$-theory $\mbox{Th}(\IR_\an)$ is denoted by $T_\an$.

\vs{0.2cm}
\noi Let $R$ be a model of $T_\an$. We call subsets and functions that are $\ma{L}_\an$-definable in $R$, again {\bf globally subanalytic}.
Note that a model $R$ of $T_\an$ is real closed and contains the reals.
Semialgebraic sets and functions are globally subanalytic.

\vs{0.2cm}
\noi Let $\ma{L}_\an^\dagger$ be the extension by definition of $\ma{L}_\an$ by a unary function symbol $^{-1}$ for multiplicative inverse and let $\ma{L}_\an^\ddagger$ be the extension by definitions of $\ma{L}_\an^\dagger$ by unary functions symbols $\sqrt[n]{}$ for $n$-th root where $n\in \IN$ with $n\geq 2$. These functions are interpreted in $\IR$ in the obvious way.
By [21], the theory $T_\an$ has quantifier elimination in the language $\ma{L}_\an^\dagger$ and is universally axiomatizable in the language $\ma{L}_\an^\ddagger$.
From this one obtains the following.

\vs{0.5cm}
{\bf 1.4 Fact}

\vs{0.1cm}
{\it Let $R$ be a model of $T_\an$ and let $A\subset R$. Then the $\ma{L}_\an^\ddagger$-substructure $\langle A\rangle_R$ of $R$ generated by $A$ is a model of $T_\an$.}

\vs{0.5cm}
\noi Examples for models of $T_\an$ are given by fields of power series.

\vs{0.5cm}
{\bf 1.5 Fact}

\vs{0.1cm}
{\it Let $\Gamma$ be an ordered abelian group that is divisible. The real closed field $\ma{R}:=\IR((t^{\Gamma}))$ has a natural expansion to a model of $T_\an$. The restricted logarithm is given by $\log_{\ma{R}}|_{]1/2,3/2[_\ma{R}}$.}

\vs{0.5cm}
\noi We write $\ma{R}_\an=\IR((t^\Gamma))_\an$ if we view the power series field over $\IR$ as a model of $T_\an$.

\vs{0.5cm}
{\bf 1.6 Proposition}
{\it \begin{itemize}
		\item[(1)]
		Let $R$ be a real closed field containing the reals. Assume that $R$ has archimedean value group. Then there is at most one $\ma{L}_\an$-structure on $R$ making it to a model of $T_\an$.
		\item[(2)] Let $R,S$ be models of $T_\an$ with archimedean value groups. Then a field embedding $\varphi:R\hookrightarrow S$ is an $\ma{L}_\an$-embedding. Assuming that $\Gamma_R\neq \{0\}$ if $\Gamma_S\neq\{0\}$ we obtain that $\varphi$ is continuous.
		\item[(3)] Let $\Gamma$ be an archimedean ordered group that is divisible and let $R$ be a real closed subfield of $\ma{R}:=\IR((t^\Gamma))$.
		Then $R$ is dense in $\langle R\rangle_\ma{R}$.
	\end{itemize}}
{\bf Proof:}

\vs{0.1cm}
	(1): Let $n\in\IN$, let $f\in\IR\langle x_1,\ldots,x_n\rangle$ and let $a=(a_1,\ldots,a_n)\in R^n$. We have to show that $\widetilde{f}(a)$ can be defined only in one way once we want to establish an $\ma{L}_\an$-structure on $R$ making it into a $T_\an$-model. By the axiomatization of $T_\an$ (see [21, Section 2]), in particular by Axiom AC4) there, it is enough to show this for $a\in\mathfrak{m}_R^n$.
	
	\noi Let $f=\sum_{\alpha\in\IN_0^n}a_\alpha x^\alpha$. For $N\in\IN$ let $f_N:=\sum_{||\alpha||\leq N}a_\alpha x^\alpha\in \IR[x_1,\ldots,x_n]$ and $g_N:=$ $f-f_N$.
	Then there is some $C_N\in\IR_{>0}$ such that $|g_N(x)|\leq C_N|x|^N$ for all $x\in [-1,1]^n$.
	We have, again by the axiomatization of $T_\an$, that
	$\widetilde{f}(a)-\widetilde{f_N}(a)=\widetilde{g_N}(a)$ for $N\in\IN$.
	Let $\gamma:=v_R(|a|)\in \Gamma_{>0}\cup\{\infty\}$ where $\Gamma:=\Gamma_R$.
	Then $v_R(\widetilde{g_N}(a))\geq N\gamma$.
	Since $\Gamma$ is archimedian we see that $\lim_{N\to \infty}\widetilde{g_N}(a)=0$ (in the order topology).
	So $\widetilde{f}(a)=\lim_{N\to\infty}\widetilde{f_N}(a)$.
	The values $\widetilde{f_N}(a)=f_N(a)$ are uniquely determined in $R$ since these functions are polynomials. So there is only one choice for $\widetilde{f}(a)$.
	
	\vs{0.2cm}
	\noi (2): Since $\varphi$ is clearly an embedding of real closed fields we have to show the following. Let $n\in\IN$ and let $f\in \IR\langle x_1,\ldots,x_n\rangle$. Then $\varphi\big(\widetilde{f}(a_1,\ldots,a_n)\big)=\widetilde{f}\big(\varphi(a_1),\ldots,\varphi(a_n)\big)$ for $a=(a_1,\ldots,a_n)\in R^n$.
	Again, it is enough to show this for $a\in\mathfrak{m}_R^n$.
	For $N\in\IN$ let $f_N$ and $g_N$ be defined as above.
	We obtain that for $N\in\IN$
	\begin{align*}
	\varphi\big(\widetilde{f}(a)\big)-\widetilde{f}\big(\varphi(a)\big)&=
	\varphi\big(\widetilde{f_N}(a)+\widetilde{g_N}(a)\big)-\Big(\widetilde{f_N}\big(\varphi(a)\big)+\widetilde{g_N}\big(\varphi(a)\big)\Big)\\
	&=\Big(\varphi\big(\widetilde{f_N}(a)\big)+\varphi\big(\widetilde{g_N}(a)\big)\Big)-\Big(\widetilde{f_N}\big(\varphi(a)\big)+\widetilde{g_N}\big(\varphi(a)\big)\Big)\\
	&=\varphi\big(\widetilde{g_N}(a)\big)-\widetilde{g_N}\big(\varphi(a)\big)=:b_N.
	\end{align*}
	Since $|a|$ is infinitesimal we obtain that $\varphi(|a|)=|\varphi(a)|$ is infinitesimal and hence \[\delta:=v_S(\varphi(|a|))\in\Gamma_S\cup\{\infty\}\] is positive. From $|\widetilde{g}_N(\varphi(a))|\leq C|\varphi(a)|^N$ we obtain that 
	\[v_S\big(\widetilde{g}_N(\varphi(a))\big)\geq N\delta.\]
	From $|\widetilde{g}_N(a)|\leq C|a|^N$ we obtain $\varphi(\widetilde{g}_N(a))\leq C\big(\varphi(|a|)\big)^N$ and therefore 
	\[v_S\big(\varphi\big(\widetilde{g_N}(a)\big)\big)\geq N\delta.\]
	Since $\Gamma_S$ is archimedean and $\delta>0$ we see that $\lim_{N\to\infty}b_N=0$.
	This shows that $\varphi\big(\widetilde{f}(a)\big)=$ $\widetilde{f}\big(\varphi(a)\big)$.
	
	\noi For the second statement we assume that $\Gamma_R\neq \{0\}$. It is enough to show that $\varphi$ is continuous at $0$. Choose an infinitesimal element $x_0$ of $R_{>0}$. Then $y_0:=\varphi(x)\in S_{>0}$ is also infinitesimal. Let $\varepsilon\in S_{>0}$.
	Since $\Gamma_S$ is archimedean we find some $n\in\IN$ such that $y_0^n<\varepsilon$.
	Setting $\delta:=x_0^n$ we obtain $|\varphi(x)|<\varepsilon$ for $x\in R$ with $|x|<\delta$.
	
	\vs{0.2cm}
	\noi (3): We define the sequence $(R_k)_{k\in\IN_0}$ recursively by $R_0:=R$ and letting
	$R_{k+1}$ be the field generated over $R_k$ by the elements $h(a)$ where $h$ is an $n$-ary function symbol in $\ma{L}_\an^\ddagger, a\in R_k^n$ and $n\in\IN$.
	In the proof of (1) we have seen that we find a sequence $(c_N)_{N\in \IN}$ in $R_k$ such that $\lim_{N\to\infty}c_N=h(a)$ if $h$ corresponds to a restricted analytic function. Dealing then with sums, products, inverses and arbitrary roots one sees that $R_k$ is dense in $R_{k+1}$ and that the value group does not enlarge.
	Since $\langle R\rangle_\ma{R}=\bigcup_{k\in\IN_0}R_k$ we obtain the claim.
	\hfill$\Box$

\vs{0.5cm}
\noi Let $\ma{L}_{\an,\exp}$ be the extension of $\ma{L}_\an$ by a unary function symbol $\exp$ for the exponential function and
let $\IR_{\an,\exp}$ be the natural $\ma{L}_{\an,\exp}$-structure on the real field. Its theory $\mathrm{Th}(\IR_{\an,\exp})$ is denoted by $T_{\an,\exp}$.
Let $\ma{L}_{\an,\exp,\log}$ be the extension by definition of $\ma{L}_{\an,\exp}$ by a unary function symbol $\log$ with the natural interpretation on $\IR$.
By [21], the theory $T_{\an,\exp}$ has quantifier elimination and is universally axiomatizable in the language $\ma{L}_{\an,\exp,\log}$.

\noi The theory $T_{\an,\exp}$ extends the theory $T_\an$.

\vs{0.2cm}
\noi An important non-standard model of $T_{\an,\exp}$ is the {\bf field  of logarithmic-exponential series} (or field of $LE$-series) $\IR((t))^{\mathrm{LE}}$ introduced by Van den Dries et al. [22].
It contains series of the form
\[t^{-1}e^{1/t}+2e^{1/t}+t^{-1/2}-\log t+6+t+2t^2+\ldots+e^{-1/t^2}-te^{-1/t^2}+\ldots+e^{-e^{1/t}}.\]
We need the following:

\vs{0.5cm}
{\bf 1.7 Fact}

\vs{0.1cm}
	{\it The power series field $\IR((t^\IR))_\an$ is an $\ma{L}_\an$-substructure of $\IR((t))^{\mathrm{LE}}$.
	The logarithm on $\IR((t))^{\mathrm{LE}}$ extends the partial logarithm on $\IR((t^\IR))$.
	For $r\in\IR$ we have that $\log(t^r)=r\log(t)$ with $\log(t)\in \IR((t))^{\mathrm{LE}}$ transcendental over $\IR((t^\IR))$.}

\subsection{Integration of parameterized definable functions}

\noi Comte, Lion and Rolin [15] (see also Lion and Rolin [37]) have shown the following seminal theorem:

\vs{0.5cm}
{\bf 1.8 Fact}

\vs{0.1cm}
{\it Let $n,k\in\IN$ with $k\leq n$ and let $q\in\IN_0$. Let $A\subset \IR^{q+n}$ be globally subanalytic such that $\dim(A_t)\leq k$ for all $t\in\IR^q$. The following holds:
	\begin{itemize}
		\item[(1)]
		The set
		\[\mathrm{Fin}_k(A):=\big\{t\in \IR^q\mid \mathrm{vol}_k(A_t)<\infty\}\]
		is globally subanalytic.
		\item[(2)]	
		There are $r\in\IN$, a real polynomial $P$ in $2r$ variables and globally subanalytic functions $\varphi_1,\ldots,\varphi_r:\mathrm{Fin}_k(A)\to \IR_{>0}$ such that
		\[\mathrm{vol}_k(A_t)=P\big(\varphi_1(t),\ldots,\varphi_r(t),\log(\varphi_1(t)),\ldots,\log(\varphi_r(t))\big)\]
		for all $t\in \mathrm{Fin}_k(A)$.
	\end{itemize}}

\vs{0.2cm}
\noi Here $\mathrm{vol}_k$ denotes the $k$-dimensional Hausdorff measure on $\IR^n$.
It is also shown that the set $\mathrm{Fin}_k(A)$ is semialgebraic if $A$ is semialgebraic.
A function is called {\bf constructible} if it is a finite sum of finite products of globally subanalytic functions and logarithms of positive globally subanalytic functions (see Cluckers and D. Miller [10, 11, 12]).
We need the following version of the above theorem (where $\lambda_n$ denotes the Lebesgue measure on $\IR^n$):

\vs{0.5cm}
{\bf 1.9 Corollary}

\vs{0.1cm}
{\it Let $n\in\IN$ and let $q\in\IN_0$.
	\begin{itemize}
		\item[(A)]
		Let $A\subset\IR^{q+n}$ be globally subanalytic. The following holds:
		\begin{itemize}
			\item[(1)]
			The set
			\[\mathrm{Fin}(A):=\big\{t\in \IR^q\mid \lambda_n(A_t)<\infty\}\]
			is globally subanalytic.
			\item[(2)]
			There is a constructible function $g:\IR^q\to \IR$ such that
			$\lambda_n(A_t)=g(t)$
			for all $t\in \mathrm{Fin}(A)$.
		\end{itemize}
		\item[(B)]
		Let $f:\IR^{q+n}\to \IR$ be globally subanalytic. The following holds:
		\begin{itemize}
			\item[(1)]
			The set
			\[\mathrm{Fin}(f):=\Big\{t\in \IR^q\;\big\vert\, \int_{\IR^n}|f_t(x)|\,d\lambda_n(x)<\infty\Big\}\]
			is globally subanalytic.
			\item[(2)]
			There is a constructible function $h:\IR^q\to\IR$ such that
			\[\int_{\IR^n}f_t(x)\,d\lambda_n(x)=h(t)\]
			for all $t\in \mathrm{Fin}(f)$.
		\end{itemize}
	\end{itemize}}
{\bf Proof:}

\vs{0.1cm}
	\noi (A) is Fact 1.8 in the case $k=n$.
	
	\vs{0.2cm}
	\noi (B) We use the following from Lebesgue integration theory. Let $g:\IR^n\to\IR_{\geq 0}$ be measurable. Then $\int_{\IR^n}g(x)\,d\lambda_n(x)=\lambda_{n+1}\big(\mathrm{subgraph}(g)\big)$.
	Applying this to $f_+$ and $f_-$ we are done by (A).
\hfill$\Box$

\vs{0.5cm}
\noi We have that $\mathrm{Fin}(f)$ is semialgebraic if $f$ is semialgebraic (see also
[32, Theorem 2.2]). Moreover in this case, one has more detailed information on the constructible functions obtained by integrating (see [33]).
But we will not need this.

\vs{0.2cm}
\noi Cluckers and D. Miller [10, 11, 12] have extended the work of Comte et al.:

\vs{0.5cm}
{\bf 1.10 Fact}

\vs{0.1cm}
{\it Let $n\in\IN$ and let $q\in\IN_0$. Let $f:\IR^{q+n}\to \IR$ be constructible. The following holds:
	\begin{itemize}
		\item[(1)] There is a constructible function $g:\IR^q\to\IR$ such that
		\[\mathrm{Fin}(f):=\Big\{t\in \IR^q\big\vert \int_{\IR^n}|f_t(x)|\,d\lambda_n(x)<\infty\Big\}\]
		equals the zero set of $g$.
		\item[(2)] There is a constructible function $h:\IR^q\to\IR$ such that
		\[\int_{\IR^n}f_t(x)\,d\lambda_n(x)=h(t)\]
		for all $t\in \mathrm{Fin}(f)$.
	\end{itemize}}

\section{Construction of the measure and the integral}

\noi Let $R$ be a real closed field with value group $\Gamma:=\Gamma_R$. We assume that $R$ contains the reals and that $\Gamma$ is {\bf archimedean}.

\subsection{Lebesgue data}

\noi
{\bf 2.1 Definition}

\vs{0.1cm}
	A {\bf Lebesgue datum} for $R$ is a tuple $\alpha=\big(s,\sigma,\tau\big)$ where $s:\Gamma\to R_{>0}$ is a section for $R$, $\sigma:R\hookrightarrow \IR((t^\Gamma))$ is a field embedding with respect to $s$ and
	$\tau:\Gamma\hookrightarrow (\IR,+)$ is an embedding of ordered groups.
	
	\noi The field embedding
	\[\Theta_\alpha:R\hookrightarrow \IR((t^\Gamma))\hookrightarrow \IR((t^\IR))\hookrightarrow \IR((t))^{\mathrm{LE}}\]
	where
	\begin{itemize}
		\item[(a)]
		$R\hookrightarrow \IR((t^\Gamma))$ is given by $\sigma$,
		\item[(b)]
		$\IR((t^\Gamma))\hookrightarrow \IR((t^\IR))$ is given by $\sum_{\gamma\in \Gamma}a_\gamma t^\gamma\mapsto \sum_{\delta\in\tau(\Gamma)} a_{\tau^{-1}(\delta)} t^{\delta}$ and
		\item[(d)]
		$\IR((t^\IR))\hookrightarrow \IR((t))^{\mathrm{LE}}$ is the inclusion
	\end{itemize}
	is called the {\bf associated Lebesgue embedding}.

\vs{0.5cm}
\noi Note that by Fact 1.1 and Fact 1.3 the field $R$ has a Lebesgue datum. In the case that $R=\IR$ there is only one Lebesgue datum
which is trivial.

\vs{0.2cm}
\noi For the rest of the section we fix a Lebesgue datum $\alpha=(s,\sigma,\tau)$ for $R$ with the associated Lebesgue embedding $\Theta:=\Theta_\alpha:R\to \IR((t))^{\mathrm{LE}}$. Via the embedding $\Theta$ we view $R$ as a subfield of $\IR((t))^{\mathrm{LE}}$. The latter field is abbreviated by $\ma{S}$ for the following.

\subsection{Construction of the measure and elementary properties}

\noi
{\bf 2.2 Construction}
	
\vs{0.1cm}
	Let $A\subset R^n$ be semialgebraic.
	We define its measure 
	\[\lambda_{R,n}(A)=\lambda_{R,n}^\alpha(A)\in \ma{S}_{\geq 0}\cup\{\infty\}\] 
	as follows.
	Take a formula $\phi(x,y)$ in the language of ordered rings, $x=(x_1,\ldots,x_n),$ $y=(y_1,\ldots,y_q)$, and a point $a\in R^q$ such that $A=\phi(R^n,a)$.
	Then the graph of the function $F:\IR^q\to \IR$ given by
	\[F(c):=\lambda_n\big(\phi(\IR^n,c)\big)\mbox{ if }\lambda_n\big(\phi(\IR^n,c)\big)<\infty\]
	and $F(c)=-1$ otherwise, is by the results mentioned in Section 1.5 defined in $\IR_{\an,\exp}$ by an $\ma{L}_{\an,\exp}$-formula $\psi(y,z)$.
	Then the formula $\psi(y,z)$ defines in $\ma{S}$ the graph of a function $F_\ma{S}:\ma{S}^q\to \ma{S}$. A routine model theoretic argument
	shows that $F_\ma{S}(a)$ does not depend on the choices of $\phi, a$ and $\psi$. This allows us to define $\lambda_{R,n}(A):=F_{\ma{S}}(a)$ if
	$F_\ma{S}(a)\geq 0$, and $\lambda_{R,n}(A)=\infty$ otherwise (that is, $F_\ma{S}(a)=-1$).

\vs{0.5cm}
\noi With a common model theoretic transfer argument we obtain the usual elementary properties of the Lebesgue measure.

\vs{0.5cm}
{\bf 2.3 Elementary properties}
{\it \begin{itemize}
		\item[(1)] {\bf Additivity:}\\
		Let $A,B\subset R^n$ be semialgebraic and disjoint. Then $\lambda_{R,n}(A\cup B)=\lambda_{R,n}(A)+\lambda_{R,n}(B)$.
		\item[(2)] {\bf Monotonicity:}\\
		Let $A,B\subset R^n$ be semialgebraic such that $A\subset B$. Then
		$\lambda_{R,n}(A)\leq \lambda_{R,n}(B)$.
		\item[(3)] {\bf Translation invariance:}\\
		Let $A\subset R^n$ be semialgebraic and let $c\in R^n$. Then $\lambda_{R,n}(A+c)=\lambda_{R,n}(A)$.
		\item[(4)] {\bf Product formula:}\\
		Let $A_1\subset R^m$ and $A_2\subset R^n$ be semialgebraic.
		Then
		\[\lambda_{R,m+n}(A_1\times A_2)=\lambda_{R,m}(A_1)\lambda_{R,n}(A_2).\]
		\item[(5)] {\bf Volume of cubes:}\\
		Let $c_j,d_j\in R$ with $c_j\leq d_j$ for $j\in\{1,\ldots,n\}$.
		Then
		\[\lambda_{R,n}\big(\prod_{j=1}^n[c_j,d_j]\big)=\prod_{j=1}^n(d_j-c_j).\]
		\item[(6)] {\bf Infinity:}\\
		We have $\lambda_{R,n}(R^n)=\infty$ for all $n\in\IN$.
	\end{itemize}}
{\bf Proof:}

\vs{0.1cm}
	For the readers' convenience we present the proof of (2):
	
	\vs{0.1cm}
	\noi We find formulas $\phi(x,y),\widetilde{\phi}(x,y)$ in the language of ordered rings, where $x=(x_1,\ldots,x_n),$ $y=(y_1,\ldots,y_q)$, and a point $a\in R^q$ such that $A=\phi(R^n,a)$ and $B=\widetilde{\phi}(\IR^n,a)$.
	Consider the functions $F:\IR^q\to \IR$ given by
	$F(c):=\lambda_n\big(\phi(\IR^n,c)\big)$ if $\lambda_n\big(\phi(\IR^n,c)\big)<\infty$
	and $F(c)=-1$ otherwise, and $\widetilde{F}:\IR^q\to \IR$ given by
	$\widetilde{F}(c):=\lambda_n\big(\widetilde{\phi}(\IR^n,c)\big)$ if $\lambda_n\big(\widetilde{\phi}(\IR^n,c)\big)<\infty$
	and $\widetilde{F}(c)=-1$ otherwise. Let $\psi(y,z)$ be an $\ma{L}_{\an,\exp}$-formula such that $\psi(\IR^{q+1})=\mathrm{graph}(F)$, and
	let $\widetilde{\psi}(y,z)$ be an $\ma{L}_{\an,\exp}$-formula such that $\widetilde{\psi}(\IR^{q+1})=\mathrm{graph}(\widetilde{F})$.
	Let $\Sigma(y)$ be the $\ma{L}_{\an,\exp}$-formula
	\[\forall t\;\forall\widetilde{t}\;\Big[\Big(\forall x\;\big(\phi(x,y)\rightarrow \widetilde{\phi}(x,y)\big)\wedge \psi(y,t)\wedge \widetilde{\psi}(y,\widetilde{t})\Big) \rightarrow \Big((t\leq \widetilde{t})\wedge (t=-1\rightarrow \widetilde{t}=-1)\Big)\Big].\]
	Then $\IR_{\an,\exp}\models \forall y\;\Sigma(y)$. So $\ma{S}\models \forall y\;\Sigma(y)$ and therefore $\Sigma(a)$ holds in $\ma{S}$. By construction of the measure, we get that
	$\lambda_{R,n}(A)\leq \lambda_{R,n}(B)$.
\hfill$\Box$

\vs{0.5cm}
\noi Property (5) gives that the measure of a cube is the naive volume. The same holds for other basic geometric subjects.

\vs{0.5cm}
{\bf 2.4 Example}

\vs{0.1cm}
{\it Let $B$ be a ball in $R^n$ with radius $r\in R_{>0}$. Then $\lambda_{R,n}(B)=\omega_n r^n$ where
	$\omega_n$ is the volume of the unit ball in $\IR^n$.
	In particular, we obtain in the case $n=2$ that $\lambda_{R,2}(B)=\pi r^2$.}

\vs{0.5cm}
\noi From semialgebraic geometry on the reals (see [32, Remark 2.1]) we obtain the following:

\vs{0.5cm}
{\bf 2.5 Proposition}
	
\vs{0.1cm}
{\it Let $A\subset R^n$ be semialgebraic. The following holds:
	\begin{itemize}
		\item[(1)]
		$\lambda_{R,n}(A)>0$ if and only if
		$\dim(A)=n$.
		\item[(2)]
		$\lambda_{R,n}(\overline{A})=\lambda_{R,n}(A)$.
	\end{itemize}}

\vs{0.2cm}
\noi A semialgebraic subset $A$ of $R^n$ is called {\bf integrable}, or is said to have {\bf finite measure} if $\lambda_{R,n}(A)<\infty$.
By $\chi_{R,n}=\chi_{R,n}^\alpha$ we denote the collection of all semialgebraic subsets of $R^n$ that are integrable.

\subsection{Construction of the integral and elementary properties}

\noi For defining integration we perform a similar construction.

\vs{0.5cm}
{\bf 2.6 Construction}

\vs{0.1cm}
	Let $f:R^n\to R_{\geq 0}$ be semialgebraic.
	We define its integral 
	\[\int_{R^n}f(x)\,dx=\int_{R^n}^\alpha f(x)\,dx\in \ma{S}_{\geq 0}\cup\{\infty\}\] as follows.
	Take a formula $\phi(x,s,y)$ in the language of ordered rings, $x=(x_1,\ldots,x_n),$  $y=(y_1,\ldots,y_q)$, and a point $a\in R^q$ such that $\mathrm{graph}(f)=\phi(R^{n+1},a)$. We choose thereby $\phi(x,s,y)$ in such a way that $\phi(\IR^{n+1},c)$ is the graph of a non-negative function
	$g_c:\IR^n\to \IR$ for every $c\in \IR^q$.
	The graph of the function $F:\IR^q\to \IR$ given by
	\[F(c):=\int_{\IR^n} g_c(x)\, dx \mbox{ if }\int_{\IR^n} g_c(x)\, dx<\infty\]	and $F(c)=-1$ otherwise, is by the results mentioned in Section 1.5 defined in $\IR_{\an,\exp}$ by an $\ma{L}_{\an,\exp}$-formula $\psi(y,z)$.
	Then the formula $\psi(y,z)$ defines in $\ma{S}$ the graph of a a function $F_\ma{S}:\ma{S}^n\to \ma{S}$. The model theoretic argument used in Construction 2.2
	shows that $F_\ma{S}(a)$ does not depend on the choices of $\phi, a$ and $\psi$. This allows us to define $\int_{R^n}f(x)\,dx:=F_{\ma{S}}(a)$ if
	$F_\ma{S}(a)\geq 0$, and $\int_{R^n}f(x)\,dx=\infty$ otherwise (that is, $F_\ma{S}(a)=-1$).
	
\vs{0.5cm}
\noi Applying the common model theoretic transfer argument we obtain the well-known connections to the Lebesgue measure.

\vs{0.5cm}
{\bf 2.7 Proposition}
{\it 	
	\begin{itemize}
		\item[(1)]
		Let $f:R^n\to R_{\geq 0}$ be semialgebraic.
		Then 
		\[\int_{R^n}f(x)\,dx=\lambda_{R,n+1}\big(\mathrm{subgraph}(f)\big).\]
		\item[(2)]
		Let $A\subset R^n$ be semialgebraic. Then 
		\[\lambda_{R,n}(A)=\int_{R^n} \mathbbm{1}_A(x)dx.\]
	\end{itemize}}

\vs{0.2cm}
\noi In view of the above we often write $\int f\,d\lambda_{R,n}$ or $\int f(x)\,d\lambda_{R,n}(x)$ for $\int_{R^n}f(x)\,dx$.

\vs{0.2cm}
\noi We extend the integral to semialgebraic functions that are not necessarily non-negative in the usual way.
Given $f:R^n\to R$ semialgebraic we have with $f_+:=\max(f,0)$ and $f_-:=\max(-f,0)$ that $f_+,f_-:R^n\to R_{\geq 0}$ are semialgebraic and that $f=f_+-f_-$.
We call $f$ {\bf integrable} if $\int f_+\,d\lambda_{R,n}<\infty$ and $\int f_-\,d\lambda_{R,n}<\infty$. We set then
\[\int f\,d\lambda_{R,n}:=\int f_+\,d\lambda_{R,n}- \int f_-\,d\lambda_{R,n}\in\ma{S}.\]
By $\ma{L}^1_{R,n}=\ma{L}_{R,n}^{1,\alpha}$ we denote the set of integrable functions.
Finally, we set
\[\mathrm{Int}_{R,n}=\mathrm{Int}^\alpha_{R,n}:\ma{L}^1_{R,n}\to \ma{S}, f\mapsto \int f\,d\lambda_{R,n}.\]

\vs{0.5cm}
{\bf 2.8 Elementary properties}

\vs{0.1cm}
{\it $\ma{L}^1_{R,n}$ is an $R$-vector space and the functional
	\[\mathrm{Int}_{R,n}:\ma{L}^1_{R,n}\to \ma{S}, f\mapsto \int f\,d\lambda_{R,n},\]
	is an $R$-linear map that is monotone.}

\vs{0.5cm}
\noi Let $A\subset R^n$ be semialgebraic and let $f:A\to R_{\geq 0}$ be semialgebraic. As usually one defines $\int_A f\,d\lambda_{R,n}$ as $\int \hat{f}d\lambda_{R,n}$ where $\hat{f}$ is the extension of $f$ by $0$ to $R^n$. Similar to above, one obtains the $R$-vector space $\ma{L}^1_{R,n}(A)=\ma{L}_{R,n}^{1,\alpha}(A)$ of semialgebraic functions integrable {\bf over} $A$ and the linear and monotone functional
$\ma{L}^1_{R,n}(A)\to \ma{S}, f\mapsto \int_A f\,d\lambda_n.$

\vs{0.2cm}
\noi From the construction of the measure and integral it is obvious that we obtain in the case $R=\IR$ the usual Lebesgue measure and integral (restricted to the semialgebraic setting).

\section{Canonicity and functoriality of the construction}

\subsection{Canonicity of the construction}

\noi
Let $R$ be a non-archimedean real closed field with archimedean value group $\Gamma:=\Gamma_R$ that contains (properly) the reals.

\noi We start by discussing which values are obtained as integrals.

\vs{0.2cm}
\noi Let $\alpha=(s,\sigma,\tau)$ be a Lebesgue datum for $R$ with the associated Lebesgue embedding \[\Theta:=\Theta_\alpha:R\to \IR((t))^{\mathrm{LE}}.\]
As above, let $\ma{R}:=\IR((t^\Gamma))$ and $\ma{S}:=\IR((t))^{\mathrm{LE}}$.
Let $\rho_\alpha$ be the embedding
\[\IR((t^\Gamma))\hookrightarrow \IR((t^\IR)), \sum_{\gamma\in \Gamma}a_\gamma t^\gamma\mapsto \sum_{\delta\in\tau(\Gamma)} a_{\tau^{-1}(\delta)} t^{\delta}.\]
We set $\ma{R}_\alpha:=\rho_\alpha(\ma{R})\subset \IR((t^\IR))$.
Note that $\ma{R}_\alpha=\ma{R}_{\alpha,\an}$ is a $T_\an$-submodel of $\IR((t^\IR))$ and that $\rho_\alpha:\ma{R}\to \ma{R}_\alpha$ is an $\ma{L}_\an$-isomorphism. Moreover, the definition of $\ma{R}_\alpha$ does only depend on the choice of the embedding $\tau:\Gamma\hookrightarrow \IR$. We have that $\Theta_\alpha(R)\subset \ma{R}_\alpha$.
We finally set
$\langle R\rangle_\alpha:=\langle\Theta_\alpha(R)\rangle_{\ma{R}_{\alpha,\an}}$. The definition of $\langle R\rangle_\alpha$ does depend on all components of $\alpha$.

\vs{0.2cm}
\noi Let $X:=\log(t^{-1})\in \ma{S}$. From the construction of the field of LE-series it follows that $X$ is transcendental over $\IR((t^\IR))$ (see Fact 1.7).
We have that $\IR<X<t^{\IR_{<0}}$.

\newpage
{\bf 3.1 Definition}

\vs{0.1cm}
	We call $\ma{R}_\alpha[X]$ the {\bf Lebesgue algebra of $R$} and
	$\langle R\rangle_\alpha[X]$ the {\bf volume algebra of $R$ with respect to $\alpha$}.

\vs{0.5cm}
\noi The Lebesgue algebra and the volume algebra are $R$-algebras by the homomorphism $\Theta_\alpha$.
The ordering on $\ma{S}$ equips the $R$-subalgebras $\ma{R}_\alpha[X]$ and $\langle R\rangle_\alpha[X]$ with an ordering. Note that these are induced by the above cut given by $X$.
Moreover, we can identify $\ma{R}_\alpha[X]$ and $\langle R\rangle_\alpha[X]$ with the polynomial ring over $\ma{R}_\alpha$ and $\langle R\rangle_\alpha$,
respectively. Note also
that $\langle R\rangle_\alpha[X]$ is an $R$-subalgebra of $\ma{R}_\alpha[X]$.

\vs{0.5cm}
{\bf 3.2 Lemma}
	{\it \begin{itemize}
		\item[(1)]
		Let $f\in \ma{R}_{>0}$. Then $\log(\rho_\alpha(f))\in -\tau\big(\mathrm{ord}(f)\big) X+\ma{O}_{\ma{R}_\alpha}$.
		\item[(2)]
		Let $f\in R_{>0}$. Then $\log(\Theta_\alpha(f))\in -\tau(v_R(f))X+\ma{O}_{\langle R\rangle_\alpha}$.
	\end{itemize}}

\vs{0.1cm}
{\bf Proof:}

\vs{0.1cm}
	(1): Let $\gamma:=\mathrm{ord}(f)\in\Gamma$ and let $g:=t^{-\gamma}f\in \ma{R}$.
	Then $\rho_\alpha(g)\in \IR_{>0}+\mathfrak{m}_{\ma{R}_\alpha}$. Let $r:=\tau(\gamma)\in\IR$. By Section 1.3 and Fact 1.7 we obtain
	\[\log(\rho_\alpha(f))=\log(t^r)+\log(\rho_\alpha(g))=-r\log(t^{-1})+\log(\rho_\alpha(g))\in -\tau\big(\mathrm{ord}(f)\big) X+\ma{O}_{\ma{R}_\alpha}.\]
	
	\vs{0.2cm}
	\noi (2): Let $f\in R_{>0}$ and let $\gamma:=v_R(f)\in\Gamma$. Since $t^\gamma\in\sigma(R)$ we have that \[h:=t^{-r}\rho_\alpha(\sigma(f))\in\rho_\alpha(\sigma(R))\subset \langle R\rangle_\alpha\]
	where $r:=\tau(\gamma)$. Since $h\in \IR_{>0}+\ma{O}_{\langle R\rangle_\alpha}$ and $\langle R\rangle_\alpha$ is a model of $T_\an$ we get that $\log(h)\in$ $ \ma{O}_{\langle R\rangle_\alpha}$. We are done since $\log(\Theta_\alpha(f))=-rX+\log(h)$.
\hfill$\Box$

\vs{0.5cm}
{\bf 3.3 Proposition}

\vs{0.1cm}
	{\it Let $A$ be a semialgebraic subset of $R^n$ that is integrable. Then 
	\[\lambda^\alpha_{R,n}(A)\in \langle R\rangle_\alpha[X]\] 
	with $\deg\big(\lambda^\alpha_{R,n}(A)\big)<n$.}

\vs{0.1cm}
{\bf Proof:}

\vs{0.1cm}
	Let $F:\IR^q\to \IR$ be a function obtained according to Construction 2.2.
	By Corollary 1.9
	there are $r\in\IN$, a polynomial $Q$ in $r$ variables over the ring of globally subanalytic functions on $\IR^q$ and globally subanalytic functions $\varphi_1,\ldots,\varphi_r$ on $\IR^q$ that are positive such that $F(y)=Q\big(\log(\varphi_1(y)),\ldots, \log(\varphi_r(y))\big)$
	for all $y\in \IR^q$.
	
	\vs{0.1cm}
	\noi {\it Claim:} One can assume that the total degree of $Q$ is bounded by $n-1$.
	
	\vs{0.05cm}
	\noi {\it Proof of the claim:}
	The case $n=1$ is an easy consequence of o-minimality (see [32, Theorem 2.3]). The higher dimensional cases follow inductively from Fubini's theorem by the inductive integration procedure in [10].
	\hfill$\Box_{\mathrm{Claim\,1}}$
	
	\vs{0.2cm}
	\noi So we assume that the total degree of $Q$ is bounded by $n-1$.
	Since $\langle R\rangle_\alpha$ is a model of $T_\an$ every globally subanalytic function $\varphi:\IR^q\to \IR$ has a canonical lifting $\varphi_{\langle R\rangle_\alpha}:\langle R\rangle_\alpha^q\to \langle R\rangle_\alpha$ (in particular,  $\varphi_{\langle R\rangle_\alpha}(a)\in\langle R\rangle_\alpha$ for every $a\in \langle R\rangle_\alpha^q$).
	The assertion follows now from Construction 2.2 and Lemma 3.2.
\hfill$\Box$

\vs{0.5cm}
\noi We discuss whether and how our construction of the semialgebraic Lebesgue measure and Lebesgue integral does depend on the choice of the
Lebesgue datum.
The definition of integrable sets and functions does not as we show below.

\vs{0.5cm}
{\bf 3.4 Proposition}

\vs{0.1cm}
	{\it Let $\alpha,\beta$ be Lebesgue data for $R$.
	The following holds:
	\begin{itemize}
		\item[(1)] Let $A$ be a semialgebraic subset of $R^n$. Then $A$ is integrable with respect to $\alpha$ if and only if it is integrable with respect to $\beta$.
		\item[(2)] Let $f:R^n\to R$ be semialgebraic. Then $f$ is integrable with respect to $\alpha$ if and only if it is integrable with respect to $\beta$.
	\end{itemize}}
{\bf Proof:}

\vs{0.1cm}
By Proposition 2.7 it suffices to show (1).
	Let $F:\IR^q\to\IR$ and $a\in R^q$ be according to Construction 2.2.
	Let $B:=F^{-1}(-1)$. Then $B$ is a semialgebraic subset of $\IR^q$ by Section 1.5.
	Since $\Theta_\alpha,\Theta_\beta:R\hookrightarrow \ma{S}$ are field embeddings over the reals we get that
	$\Theta_\alpha(a)\in B_\ma{S}$ if and only if $\Theta_\beta(a)\in B_\ma{S}$.
	This gives (1).
\hfill$\Box$

\vs{0.5cm}
\noi From now on we write $\chi_{R,n}$ and $\ma{L}^1_{R,n}$ for the set of integrable semialgebraic subsets of $R^n$ and of integrable semialgebraic functions
$R^n\to R$, respectively. For the following we restrict a semialgebraic measure $\lambda_{R,n}^\alpha$ to $\chi_{R,n}$.

\vs{0.2cm}
\noi We introduce a natural notion of equivalence between the constructions obtained from choosing different Lebesgue data.

\vs{0.5cm}
{\bf 3.5 Definition}
\begin{itemize}
		\item[(a)]
		Let $\alpha,\beta$ be two Lebesgue data for $R$.
		We say that the semialgebraic Lebesgue measure with respect to $\alpha$ and the semialgebraic Lebesgue measure with respect to $\beta$ are
		{\bf isomorphic} if there is a ring isomorphism $\Phi:\ma{R}_\alpha[X]\stackrel{\cong}{\longrightarrow} \ma{R}_\beta[X]$ of the Lebesgue algebras such that $\lambda^\beta_{R,n}=\Phi\circ\lambda_{R,n}^\alpha$ for every $n\in\IN$.
		We call such an isomorphism a {\bf Lebesgue isomorphism} between the Lebesgue data $\alpha$ and $\beta$.
		\item[(b)]
		We say that the semialgebraic Lebesgue measure for $R$ is
		{\bf unique up to isomorphism} if the semialgebraic Lebesgue measures with respect to any Lebesgue data for $R$ are isomorphic.
	\end{itemize}

\vs{0.2cm}
{\bf 3.6 Remark}

\vs{0.1cm}
	{\it A Lebesgue isomorphism $\Phi:\ma{R}_\alpha[X]\stackrel{\cong}{\longrightarrow}\ma{R}_\beta[X]$ is an $R$-algebra isomorphism which is order preserving.}

\vs{0.1cm}
{\bf Proof:}

\vs{0.1cm}
	Let $\alpha=(s_\alpha,\sigma_\alpha,\tau_\alpha)$ and $\beta=(s_\beta,\sigma_\beta,\tau_\beta)$.
	
	\vs{0.2cm}
	\noi a) Let $a\in R_{>0}$. We obtain by 2.3(5) and Construction 2.2 that
	\[\Phi\big(\Theta_\alpha(a)\big)=\Phi\big(\lambda_{R,1}^\alpha([0,a])\big)
	=\lambda_{R,1}^\beta([0,a])=\Theta_\beta(a).\]
	This shows that $\Phi$ is an $R$-algebra isomorphism.
	
	\vs{0.2cm}
	\noi b) Looking at the invertible elements of the polynomial ring we get that $\Phi(\ma{R}_\alpha)=\ma{R}_\beta$. Since these fields are real closed we obtain that
	$\Phi|_{\ma{R}_\alpha}$ is order preserving.
	Let $\gamma\in \Gamma_{<0}$ and let $c\in R_{>0}$ with $v_R(c)=\gamma$. Let
	\[A:=\big\{(x,y)\in R^2\mid 1\leq x\leq c, 0\leq y\leq 1/x\big\}.\]
	By Proposition 2.7(1), Construction 2.6 and Lemma 3.2(2) we obtain that
	\[\lambda_{R,2}^\alpha(A)=\int_{[1,c]}^\alpha dx/x=\log\big(\Theta_\alpha(c)\big)
	= -\tau_\alpha(\gamma)X+h_\alpha\]
	and that
	\[\lambda_{R,2}^\beta(A)=\int_{[1,c]}^\beta dx/x=\log\big(\Theta_\beta(c)\big)
	= -\tau_\beta(\gamma)X+h_\beta\]
	where $h_\alpha\in \ma{O}_{\langle R\rangle_\alpha}$ and $h_\beta\in \ma{O}_{\langle R\rangle_\beta}$.
	This shows that there are $r\in \IR_{>0}$ and $g\in \ma{O}_{\ma{R}_\beta}$ such that $\Phi(X)=rX+g$. By the ordering given on $\ma{R}_\alpha[X]$ respectively $\ma{R}_\beta[X]$ we conclude that $\Phi$ is order preserving.
\hfill$\Box$

\vs{0.5cm}
{\bf 3.7 Proposition}

\vs{0.1cm}
	{\it Assume that two Lebesgue data $\alpha$ and $\beta$ differ only by the embedding of the value group $\Gamma$ into the group of reals.
	Then the semialgebraic Lebesgue measure with respect to $\alpha$ is isomorphic to the semialgebraic Lebesgue measure with respect to $\beta$.}

\newpage
{\bf Proof:}

\vs{0.1cm}
	Let $\alpha=(s,\sigma,\tau_\alpha)$ and $\beta=(s,\sigma,\tau_\beta)$.
	We have to find an isomorphism
	\[\Phi:\ma{R}_\alpha[X]\to \ma{R}_\beta[X]\] such that $\Phi\big(F_\ma{S}(\Theta_\alpha(a))\big)=F_\ma{S}\big(\Theta_\beta(a)\big)$ for a function  $F:\IR^q\to \IR$ and a tuple $a\in R^q$ obtained by applying Construction 2.2 to a semialgebraic subset of some $R^n$. By Corollary 1.9(A) we know that such an $F$ is constructible.
	Via $\sigma$ we can identify $R$ with a subfield of $\ma{R}$.
	By [40, Satz 5 in I \S 3] we find some $r\in\IR_{>0}$ such that that $\tau_\beta=r\tau_\alpha$.
	
	\vs{0.2cm}
	\noi {\it Claim 1:}
	The field automorphism
	\[G:\IR((t^\IR))\to\IR((t^\IR)), \sum_{\alpha\in\IR}a_{\alpha}t^{\alpha}\mapsto \sum_{\alpha\in\IR}a_{\alpha}t^{r\alpha}\]
	is an automorphism of the $\ma{L}_\an$-structure $\IR((t^\IR))$.
	
	\vs{0.1cm}
	\noi {\it Proof of Claim 1:}
	This follows from the definition of the $\ma{L}_\an$-structure on power series fields (see [21, Section 2]) or from Proposition 1.6(2).
	\hfill$\Box_{\mathrm{Claim\, 1}}$
	
	\vs{0.2cm}
	\noi We have that $\rho_\beta=G\circ \rho_\alpha$. So $G$ restricts to an isomorphism $H:\ma{R}_{\alpha,\an}\to \ma{R}_{\beta,\an}$
	with $\rho_\beta=H\circ \rho_\alpha$.
	The map
	$\Phi:\ma{R}_\alpha[X]\to \ma{R}_\beta[X]$ extending $H$ and sending $X$ to $rX$ (note that $\log(t^{-r})=r\log(t^{-1}$)) is an isomorphism of the Lebesgue algebras. By the definition of a constructible function we are done by Claim 1 and establishing the following Claim 2:
	
	\vs{0.2cm}
	\noi {\it Claim 2:}
	Let $x\in \ma{R}_{\alpha}$ be positive. Then $\log\big(\Phi(x)\big)=\Phi\big(\log(x)\big)$.
	
	\vs{0.1cm}
	\noi {\it Proof of Claim 2:}
	Let $\delta\in\IR,a\in \IR_{>0}$ and $g\in\mathfrak{m}_{\ma{R}_\alpha}$ such that $x=at^\delta(1+g)$.
	Since $H(g)\in\mathfrak{m}_{\ma{R}_\beta}$ we obtain
	\begin{align*}
	\log\big(\Phi(x)\big)&=
	\log\big(H(at^\delta(1+g))\big)\\
	&=\log\big(at^{r\delta}(1+H(g))\big)\\
	&=\log(a)-r\delta X+L(H(g)))\\
	&\stackrel{\mathrm{Claim\, 1}}{=}
	\log(a)-r\delta X+ H(L(g))\\
	&=\Phi\big(\log(a)-\delta X+ L(g)\big)\\
	&=\Phi\big(\log(x)\big).
	\end{align*}
	where $L$ denotes the logarithmic series. This shows Claim 2.
	\hfill$\Box_{\mathrm{Claim\, 2}}$
	
	\hfill$\Box$
	
\vs{0.5cm}
{\bf 3.8 Proposition}

\vs{0.1cm}
{\it 
	Let $s,s'$ be sections for $R$ and let $\sigma,\sigma':R\to \ma{R}$ be embeddings with respect to $s$ and $s'$, respectively.
	Then there is a valuation preserving $\ma{L}_\an$-automorphism $K:\ma{R}_\an\to \ma{R}_\an$ such that $\sigma'=K\circ\sigma$.}

\vs{0.1cm}
{\bf Proof:}

\vs{0.1cm}
We have that $\ma{R}$ is the immediate maximal extension of $R$ (see Kaplansky [34] and [40, III \S 3]). Since $\sigma,\sigma'$ are valuation preserving by Fact 1.3 we get by [34, Theorem 5] that there is a valuation preserving field automorphism $K:\ma{R}\to \ma{R}$ such that $\sigma'=K\circ\sigma$. By Proposition 1.6(2) we obtain that $K$ is an $\ma{L}_\an$-automorphism of $\ma{R}_\an$.
\hfill$\Box$

\vs{0.5cm}
\noi In general, the construction of the Lebesgue measure does only depend on the choice of a section for the given real closed field.

\vs{0.5cm}
{\bf 3.9 Theorem}

\vs{0.1cm}
{\it 
	Let $\alpha,\beta$ be Lebesgue data having the same section. Then the semialgebraic Lebesgue measure with respect to $\alpha$ is isomorphic to the semialgebraic Lebesgue measure with respect to $\beta$.}

\newpage
{\bf Proof:}

\vs{0.1cm}
Let $\alpha=(s,\sigma_\alpha,\tau_\alpha)$ and let $\beta=(s,\sigma_\beta,\tau_\beta)$.
	By Proposition 3.7 we can assume that $\tau_\alpha=\tau_\beta=:\tau$.
	Via $\tau$ we can identify $\Gamma$ with a subgroup of $(\IR,+)$ and obtain that $\ma{R}=\ma{R}_\alpha=\ma{R}_\beta$.
	By Proposition 3.8 there is a valuation preserving $\ma{L}_\an$-automorphism $K:\ma{R}_\an\to \ma{R}_\an$ such that $\sigma_\beta=K\circ\sigma_\alpha$.
	
	\vs{0.2cm}
	\noi {\it Claim 1:} We have that $K(t^\gamma)=t^\gamma$ for all $\gamma\in\Gamma$.
	
	\vs{0.1cm}
	\noi {\it Proof of Claim 1:}
	Since $\sigma_\alpha$ and $\sigma_\beta$ are embeddings with respect to $s$ we have that
	\[K(t^\gamma)=K\big(\sigma_\alpha(s(\gamma))\big)=\sigma_\beta(s(\gamma))=t^\gamma\]
	for all $\gamma\in\Gamma$.
	\hfill$\Box_{\mathrm{Claim\, 1}}$
	
	\vs{0.1cm}
	\noi Mapping $X\mapsto X$ we extend $K$ to an isomorphism $\Phi:\ma{R}[X]\to \ma{R}[X]$.
	As in the proof of Proposition 3.7 we are done once the following Claim 2 is established:
	
	\vs{0.2cm}
	\noi {\it Claim 2:}
	Let $x\in \ma{R}_{>0}$. Then $\log\big(\Phi(x)\big)=\Phi\big(\log(x)\big)$.
	
	\vs{0.1cm}
	\noi {\it Proof of Claim 2:}
	Let $\gamma\in\Gamma,a\in \IR_{>0}$ and $h\in\mathfrak{m}_{\ma{R}}$ such that $x=at^\gamma(1+h)$.
	We obtain
	\begin{align*}
	\log\big(\Phi(x)\big)&=
	\log\Big(K\big((at^\gamma(1+h)\big)\Big)\\
	&=\log\big(K(t^\gamma)\big)+\log(a)+L\big(K(h)\big)\\
	&\stackrel{\mathrm{Claim\,1}}{=}
	\log\big(t^\gamma\big)+\log(a)+L\big(K(h)\big)\\
	&=-\gamma X+\log(a)+K\big(L(h)\big)\\
	&=\Phi\big(\log(a)-\gamma X+ L(h)\big)\\
	&=\Phi\big(\log(x)\big).
	\end{align*}
	where $L$ denotes the logarithmic series.
	\hfill$\Box_{\mathrm{Claim\, 2}}$
	
	\hfill$\Box$

\vs{0.5cm}
\noi Theorem 3.9 is the best we can hope for in the case that the value group is not isomorphic to the rationals.

\vs{0.5cm}
{\bf 3.10 Theorem}

\vs{0.1cm}
{\it 
	Assume that $\mathrm{rank}(\Gamma)>1$.
	Then the semialgebraic Lebesgue measure is not unique up to isomorphism.}

\vs{0.1cm}
{\bf Proof:}

\vs{0.1cm}
	Let $\gamma\in \Gamma_{<0}$ and let $\tau:\Gamma\hookrightarrow \IR$ be an embedding
	with $\tau(\gamma)=-1$.
	Choose an element $\delta\in\Gamma_{<0}$ that is linearly independent from $\gamma$ and let $\zeta:=-\tau(\delta)$. Then $\zeta$ is irrational.
	We choose distinct elements $a,b,b'\in R_{>0}$ such that $v_R(a)=\gamma$ and $v_R(b)=v_R(b')=\delta$.
	Let $s$ be a section for $R$ that maps $\gamma$ to $a$ and $\delta$ to $b$.
	Let $s'$ be a section for $R$ that maps $\gamma$ to $a$ and $\delta$ to $b'$.
	This is possible since $\gamma$ and $\delta$ are linearly independent.
	Let $\sigma, \sigma':R\to\ma{R}$ be field embeddings with respect to $s$ and $s'$, respectively.
	Then the semialgebraic measure with respect to the Lebesgue datum  $\alpha=(s,\sigma,\tau)$ and the semialgebraic measure with respect to the Lebesgue datum $\beta=(s',\sigma',\tau)$ are not isomorphic. To show this
	let
	\[A:=\big\{(x,y)\in R^2\mid 1\leq x\leq a, 0\leq y\leq 1/x\big\}\]
	and
	\[B:=\big\{(x,y)\in R^2\mid 1\leq x\leq b, 0\leq y\leq 1/x\big\}.\]
	We have that $\Theta_\alpha(a)=\Theta_\beta(a)=t^{-1}$.
	We also have that $\Theta_\alpha(b)=\Theta_\beta(b')=t^{-\zeta}$. Let $h:=b/b'$. Then $h$ is a unit in the valuation ring of $\ma{O}_R$; i.e., $v_R(h)=0$. Since $b\neq b'$ we have that $h\neq 1$. We obtain
	\[\Theta_\beta(b)=\Theta_\beta(b'h)=t^{-\zeta}\Theta_\beta(h)\]
	where $\mathrm{ord}\big(\Theta_\beta(h)\big)=0$.
	Computing the Lebesgue measures of $A$ and $B$ with respect to $\alpha$ we obtain
	\[\lambda_{R,2}^\alpha(A)=\int_{[1,a]}^\alpha dx/x=\log\big(\Theta_\alpha(a)\big)=\log(t^{-1})=X\]
	and
	\[\lambda_{R,2}^\alpha(B)=\int_{[1,b]}^\alpha dx/x=\log\big(\Theta_\alpha(b)\big)=\log(t^{-\zeta})=\zeta X.\]
	Doing the computation with respect to $\beta$ we get
	\[\lambda_{R,2}^\beta(A)=\int_{[1,a]}^\beta dx/x=\log\big(\Theta_\beta(a)\big)=\log(t^{-1})=X\]
	and
	\[\lambda_{R,2}^\beta(B)=\int_{[1,b]}^\beta dx/x=\log\big(\Theta_\beta(b)\big)=\log\big(t^{-\zeta}\Theta_\beta(h)\big)=\zeta X+g\]
	where $g:=\log(\Theta_\beta(h))\neq 0$.
	Hence there can be no homomorphism $\Phi:\ma{R}_\alpha[X]\to \ma{R}_\beta[X]$ with $\lambda_{R,2}^\beta=\Phi\circ\lambda_{R,2}^\alpha$.
\hfill$\Box$

\vs{0.5cm}
\noi In the case that the value group is isomorphic to the group of rationals we obtain the best possible result.

\vs{0.5cm}
{\bf 3.11 Theorem}	

\vs{0.1cm}
{\it 	
	Assume that $\Gamma\cong \IQ$. Then the semialgebraic Lebesgue measure is unique up to isomorphism.}

\vs{0.1cm}
{\bf Proof:}

\vs{0.1cm}
	Let $\alpha=(s_\alpha,\sigma_\alpha,\tau_\alpha)$ and $\beta=(s_\beta,\sigma_\beta,\tau_\beta)$ be two Lebesgue data for $R$.
	By Proposition 3.7 we may assume that $\tau_\alpha=\tau_\beta=:\tau$.
	Via $\tau$ we can identify $\Gamma$ with a subgroup of $(\IR,+)$ and obtain that $\ma{R}=\ma{R}_\alpha=\ma{R}_\beta$.
	By Proposition 3.8 there is a valuation preserving $\ma{L}_\an$-automorphism $K:\ma{R}_\an\to \ma{R}_\an$ such that $\sigma_\beta=K\circ\sigma_\alpha$.
	Since $K$ is an automorphism of the real field $\ma{R}$ we have that $K(x^q)=\big(K(x)\big)^q$ for all $x\in\ma{R}_{>0}$ and all $q\in\IQ$.
	
	\noi Since $\Gamma$ is isomorphic to $\IQ$ there is an $r\in \IR_{>0}$ such that $\Gamma=r\IQ$.
	Note that the choice of $r$ depends only on its class in $\IR/\IQ$.
	We have that $-r\in \Gamma$. Since $t^{-r}\in \sigma_\alpha(R)$ we find a (uniquely determined) $x^*\in R$ such that $\sigma_\alpha(x^*)=t^{-r}$. Then $\mathrm{ord}(\sigma_\beta(x^*))=-r$. Since $t^{-r}\in \sigma_\beta(R)$ we find $a^*\in\IR_{>0}$ and $h^*\in \mathfrak{m}_R$ such that $\sigma_\beta(x^*)=a^* t^{-r}(1+\sigma_\beta(h^*))$. Then $f^*:=\log\big(a^*(1+\sigma_\beta(h^*))\big)\in \ma{O}_{\langle R\rangle_\beta}$ by Lemma 3.2(2).
	Mapping $X\mapsto X+f^*/r$ we extend $K$ to an isomorphism $\Phi:\ma{R}[X]\to \ma{R}[X]$.
	As above we are done by establishing the following Claim.
	
	\vs{0.1cm}
	\noi {\it Claim:}
	Let $x\in \ma{R}_{>0}$. Then $\log\big(\Phi(x)\big)=\Phi\big(\log(x)\big)$.
	
	\vs{0.05cm}
	\noi {\it Proof of the claim:}
	Let $\gamma\in\Gamma,a\in \IR_{>0}$ and $h\in\mathfrak{m}_{\ma{R}}$ such that $x=at^\gamma(1+h)$.
	Let $q_\gamma:=\gamma/r\in\IQ$.
	We obtain
	\begin{align*}
	\log\big(\Phi(x)\big)&=
	\log\Big(K\big((t^{-r})^{-q_\gamma}a(1+h)\big)\Big)\\
	&=
	\log\Big(\big(K(t^{-r})\big)^{-q_\gamma}K\big(a(1+h)\big)\Big)\\
	&=-q_\gamma\log\big(\sigma_\beta(x^*)\big)+\log\Big(K\big(a(1+h)\big)\Big)\\
	&=-q_\gamma(rX+f^*)+\log\Big(a\big(1+K(h)\big)\Big)\\
	&=-\gamma(X+f^*/r)+\log(a)+L\big(K(h)\big)\\
	&=
	-\gamma(X+f^*/r)+\log(a)+K\big(L(h)\big)\\
	&=\Phi\big(\log(a)-\gamma X+ L(h)\big)\\
	&=\Phi\big(\log(x)\big).
	\end{align*}
	where $L$ denotes the logarithmic series.
	\hfill$\Box_{\mathrm{Claim}}$
	
\hfill$\Box$

\vs{0.5cm}
\noi One could also define a {\bf volume isomorphism} $\langle R\rangle_\alpha[X]\stackrel{\cong}{\longrightarrow}\langle R\rangle_\beta[X]$ between two Lebesgue data $\alpha$ and $\beta$. But this does not lead to
a new notion of isomorphism:

\newpage
{\bf 3.12 Proposition}

\vs{0.1cm}
{\it 
	Let $\alpha,\beta$ be Lebesgue data for $R$. The following holds:
	\begin{itemize}
		\item[(1)] A Lebesgue isomorphism $\ma{R}_\alpha[X]\stackrel{\cong}{\longrightarrow}\ma{R}_\beta[X]$ restricts to a volume isomorphism\\ $\langle R\rangle_\alpha[X]\stackrel{\cong}{\longrightarrow}\langle R\rangle_\beta[X].$
		\item[(2)]  A volume isomorphism $\langle R\rangle_\alpha[X]\stackrel{\cong}{\longrightarrow}\langle R\rangle_\beta[X]$ can be extended to a Lebesgue isomorphism $\ma{R}_\alpha[X]\stackrel{\cong}{\longrightarrow}\ma{R}_\beta[X]$.
	\end{itemize}}
{\bf Proof:}

\vs{0.1cm}
	Let $\alpha=(s_\alpha,\sigma_\alpha,\tau_\alpha)$ and $\beta=(s_\beta,\sigma_\beta,\tau_\beta)$.
	By [40, Satz 5 in I \S 5] we find some $r\in \IR_{>0}$ such that
	$\tau_\beta=r\tau_\alpha$. Let $\eta:\tau_\alpha(\Gamma)\to \tau_\beta(\Gamma), \delta\mapsto r\delta$.
	
	\vs{0.2cm}
	\noi (1): We denote the given Lebesgue isomorphism by $\Phi$.
	We have that $\Phi(\ma{R}_\alpha)=\ma{R}_\beta$.
	Let $\varphi:=\Phi|_{\ma{R}_\alpha}$.
	By Remark 3.6 we have that
	$\varphi\big(\Theta_\alpha(R)\big)\subset \Theta_\beta(R)$.
	By Proposition 1.6(2) we know that $\varphi$ is an $\ma{L}_\an$-isomorphism. Hence
	$\varphi\big(\langle R\rangle_\alpha\big)=\big\langle \Phi\big(\Theta_\alpha(R)\big)\big\rangle_{\ma{R}_\beta}$. So we obtain
	$\varphi\big(\langle R\rangle_\alpha\big)\subset \langle R\rangle_\beta$.
	Using this and the observation b) in the proof of Remark 3.6 we get that there are $\widetilde{r}\in\IR_{>0}$ and $g\in \ma{O}_{\langle R\rangle_\beta}$ such that $\Phi(X)=\widetilde{r}X+g$.
	So $\Phi(X)\in \langle R\rangle_\beta[X]$.
	We obtain that $\Phi\big(\langle R\rangle_\alpha[X]\big)\subset\langle R\rangle_\beta[X]$ and by symmetry we obtain equality.
	
	\vs{0.2cm}
	\noi (2): We denote the given volume isomorphism by $\Psi$. We have that $\Psi(\langle R\rangle_\alpha)=\langle R\rangle_\beta$ and see similar to the proof of Proposition 3.7 that $\Psi(X)=rX$. Let $\psi:=\Psi|_{\langle R\rangle_\alpha}$.
	
	\vs{0.2cm}
	\noi {\it Special case:} $\sigma_\alpha=\sigma_\beta$.
	
	\vs{0.1cm}
	\noi Then also $s_\alpha=s_\beta$. Via $\sigma:=\sigma_\alpha$ we can identify $R$ with a subfield of $\ma{R}$.
	Let
	\[G:\IR((t^\IR))\to\IR((t^\IR)), \sum_{\alpha\in\IR}a_\alpha t^\alpha\mapsto \sum_{\alpha\in\IR}a_\alpha t^{r\alpha}.\]
	By Remark 3.6 we have that $\psi|_R=G|_R$. By Proposition 1.6(3) we have that $R$ is dense in $\langle R\rangle_\alpha$ and by Proposition 1.6(2) that $\psi$ and $G$ are continuous.
	Hence we obtain that $\psi=G|_{\langle R\rangle_\alpha}$.
	The isomorphism $\Phi:\ma{R}_\alpha[X]\to \ma{R}_\beta[X]$ with $\Phi|_{\ma{R}_\alpha}=G|_{\ma{R}_\alpha}$ and sending $X\to rX=\Psi(X)$ is clearly a Lebesgue isomorphism between $\alpha$ and $\beta$ and extends $\Psi$.
	
	\vs{0.2cm}
	\noi {\it General case:}
	By the special case we can assume that $\tau_\alpha=\tau_\beta$. Therefore we identify $\Gamma$ with a subgroup of $\IR$ and obtain $\ma{R}=\ma{R}_\alpha=\ma{R}_\beta$. The maps $\Theta_\alpha$ and $\Theta_\beta$ are valuation preserving. Since $\psi\circ\Theta_\alpha=\Theta_\beta$ by Remark 3.6 we get that $\psi|_{\Theta_\alpha(R)}$ is valuation preserving.
	Since $\Theta_\alpha(R)$ and $\langle R\rangle_\alpha$ have the same value group we see that $\psi$ is valuation preserving.
	Since $\ma{R}$ is the immediate maximal extension of $\langle R\rangle_\alpha$ respectively $\langle R\rangle_\beta$ we get by [34, Theorem 5] that there is a valuation preserving field automorphism $K:\ma{R}\to\ma{R}$ extending $\psi$. We finish as in the special case.
\hfill$\Box$

\vs{0.5cm}
\noi To obtain a kind of weak uniqueness we introduce the following notations.

\vs{0.5cm}
{\bf 3.13 Definition}
	\begin{itemize}
		\item[(a)] Let $\alpha$ be a Lebesgue datum for $R$. Then the ordered group
		\[\ma{R}_\alpha[X]/\ma{O}_{\ma{R}_\alpha}\]  is called the {\bf reduced Lebesgue group} of $R$ with respect to $\alpha$. For $n\in\IN$, we set $\overline{\lambda}^\alpha_{R,n}:=\pi\circ\lambda_{R,n}^\alpha$ where $\pi:\ma{R}_\alpha[X]\to \ma{R}_\alpha[X]/\ma{O}_{\ma{R}_\alpha}$ is the canonical group epimorphism.
		\item[(b)]
		Let $\alpha,\beta$ be two Lebesgue data for $R$.
		We say that the reduced semialgebraic Lebesgue measure with respect to $\alpha$ and the reduced semialgebraic Lebesgue measure with respect to $\beta$ are
		{\bf isomorphic} if there is a group isomorphism \[\Psi:\ma{R}_\alpha[X]/\ma{O}_{\ma{R}_\alpha}\stackrel{\cong}{\longrightarrow} \ma{R}_\beta[X]/\ma{O}_{\ma{R}_\beta}\] of the reduced Lebesgue groups such that $\overline{\lambda}^\beta_{R,n}=\Psi\circ\overline{\lambda}_{R,n}^\alpha$ for every $n\in\IN$.
		We call such an isomorphism a {\bf reduced Lebesgue isomorphism} between the Lebesgue data $\alpha$ and $\beta$.
		\item[(c)]
		We say that the reduced semialgebraic Lebesgue measure for $R$ is
		{\bf unique up to isomorphism} if the reduced semialgebraic Lebesgue measures with respect to any Lebesgue data for $R$ are isomorphic.
	\end{itemize}

\vs{0.2cm}
{\bf 3.14 Theorem}

\vs{0.1cm}
{\it 
	The reduced semialgebraic Lebesgue measure is unique up to isomorphism.}

\vs{0.1cm}
{\bf Proof:}

\vs{0.1cm}
	Let $\alpha=(s_\alpha,\sigma_\alpha,\tau_\alpha)$ and $\beta=(s_\beta,\sigma_\beta,\tau_\beta)$ be two Lebesgue data for $R$.
	As above, we may assume that $\tau_\alpha=\tau_\beta=:\tau$ and identify $\Gamma$ via $\tau$ with a subgroup of $(\IR,+)$, obtaining thereby that  $\ma{R}=\ma{R}_\alpha=\ma{R}_\beta$.
	By Proposition 3.8 there is a valuation preserving $\ma{L}_\an$-automorphism $K:\ma{R}_\an\to \ma{R}_\an$ such that $\sigma_\beta=K\circ\sigma_\alpha$.
	Mapping $X\mapsto X$ we extend $K$ to an algebra isomorphism $\Phi:\ma{R}[X]\to \ma{R}[X]$.
	Since $\Phi(\ma{O}_{\ma{R}})=\ma{O}_{\ma{R}}$ we get that $\Phi$ induces a group isomorphism $\Psi:\ma{R}[X]/\ma{O}_{\ma{R}}\to \ma{R}[X]/\ma{O}_{\ma{R}}$.
	Since $K$ is an $\ma{L}_\an$-automorphism and since $\Phi$ is a ring isomorphism we are done, as above, once the following claim is established where $\pi:\ma{R}[X]\to \ma{R}[X]/\ma{O}_{\ma{R}}$ denotes the canonical epimorphism.
	
	\vs{0.2cm}
	\noi {\it Claim:}
	Let $x\in \ma{R}_{>0}$. Then $\pi\Big(\log\big(K(x)\big)\Big)=\Psi\Big(\pi\big(\log(x)\big)\Big)$.
	
	\vs{0.1cm}
	\noi {\it Proof of the claim:}
	Let $\gamma\in\Gamma,a\in \IR_{>0}$ and $h\in\mathfrak{m}_{\ma{R}}$ such that $x=at^\gamma(1+h)$. Since $K$ is valuation preserving we find $b\in \IR_{>0}$ and $g\in\mathfrak{m}_{\ma{R}}$ such that $K(t^\gamma)=bt^\gamma(1+g)$.
	We obtain
	\begin{align*}
	\pi\Big(\log\big(K(x)\big)\Big)&=
	\pi\Big(\log\big(bt^\gamma(1+g)\big)+\log(a)+L\big(K(h)\big)\Big)\\
	&=\pi\Big(\gamma X+\log(ab)+L(g)+L\big(K(h)\big)\Big)\\
	&=\gamma X=\Psi(\gamma X)\\
	&=\Psi\Big(\pi\big(\gamma X+\log(a)+L(h)\big)\Big)\\
	&=\Psi\Big(\pi\big(\log(x)\big)\Big).\\
	\end{align*}
	where $L$ denotes the logarithmic series.
	\hfill$\Box_{\mathrm{Claim}}$
	
\hfill$\Box$

\vs{0.5cm}
\noi The above statements hold analogously for the semialgebraic integral.

\subsection{The Lebesgue algebra and the volume algebra}

\noi As in Section 3.1 let $R$ be a non-archimedean real closed field with archimedean value group $\Gamma:=\Gamma_R$ that contains (properly) the reals.
As above we set $\ma{R}:=\IR((t^\Gamma))$.
We have seen in Theorem 3.9 that the construction of the Lebesgue measure and the Lebesgue integral depends only on the choice of the section, not on the choice of the embedding with respect to the section and not on the choice of embedding of the value group into the reals.
Since the definition of $\ma{R}_\alpha$ depends only on the latter the following definitions are
justified:

\vs{0.5cm}
{\bf 3.15 Definition}
	\begin{itemize}
		\item[(1)]
		We call the polynomial algebra $\ma{R}[X]$ over $\ma{R}$ in one variable the {\bf Lebesgue algebra} of $R$.
		\item[(2)]
		Let $s$ be a section for $R$. For
		$n\in\IN$ we call
		\[\lambda_{R,n}=\lambda_{R,n}^s:\big\{\mbox{semialgebraic subsets of }R^n\big\}\to \ma{R}[X]\cup\{\infty\},A\mapsto \lambda_{R,n}(A),\]
		and
		\[\mathrm{Int}_{R,n}=\mathrm{Int}^s_{R,n}:\ma{L}^1_{R,n}\to \ma{R}[X],f\mapsto \int_{R^n}f\,d\lambda_{R,n},\]
		the {\bf semialgebraic measure and integral on $R^n$ with respect to $s$}.
	\end{itemize}

\vs{0.2cm}
{\bf 3.16 Remark}
	{\it \begin{itemize}
		\item[(1)] Let $A$ be a semialgebraic subset of $R^n$ with finite measure. Then \[\mathrm{deg}\big(\lambda_{R,n}(A)\big)<n.\]
		\item[(2)] Let $f:R^n\to R$ be a semialgebraic function that is integrable. Then 
		\[\mathrm{deg}\Big(\int_{R^n}f\,d\lambda_{R,n}\Big)\leq n.\]
	\end{itemize}
	Moreover, the degree does not depend on the choice of the section.}

\vs{0.1cm}
{\bf Proof:}

\vs{0.1cm}
	(1) has been shown in Proposition 3.3. (2) follows then from Proposition 2.7(1).
	That the degree does not depend on the choice of the section follows from Constructions 2.2 and 2.6.
	\hfill$\Box$

\vs{0.5cm}
\noi If the value group is isomorphic to the rationals we have seen in Theorem 3.10 that the semialgebraic Lebesgue measure and integral is unique up to isomorphisms. Embedding $R$ into $\ma{R}$ and setting $\langle R\rangle:=\langle R\rangle_\ma{R}$, the following definitions are
justified:

\vs{0.5cm}
{\bf 3.17 Definition}

\vs{0.1cm}
	Assume that $\Gamma\cong \IQ$.
	\begin{itemize}
		\item[(1)]
		We call the polynomial algebra $\langle R\rangle[X]$ over $\langle R\rangle$ in one variable the {\bf volume algebra} of $R$.
		\item[(2)]
		For $n\in\IN$ we call
		\[\lambda_{R,n}:\big\{\mbox{semialgebraic subsets of }R^n\big\}\to \langle R\rangle[X]\cup\{\infty\},A\mapsto \lambda_{R,n}(A),\]
		and
		\[\mathrm{Int}_{R,n}:\ma{L}^1_{R,n}\to \langle R\rangle[X],f\mapsto \int_{R^n}f\,d\lambda_{R,n},\]
		{\bf the  semialgebraic measure and integral on $R^n$}.
	\end{itemize}

\vs{0.2cm}
{\bf 3.18 Theorem}

\vs{0.1cm}
	{\it Assume that $\Gamma\cong \IQ$ and that $R$ can be made to a model of $T_\an$. Then the volume algebra of $R$ is the polynomial algebra $R[X]$ over $R$.}

\vs{0.1cm}
{\bf Proof:}

\vs{0.1cm}
	By Proposition 1.6(1) we see that $\langle R\rangle=R$.
\hfill$\Box$

\vs{0.5cm}
{\bf 3.19 Main example}

\vs{0.1cm}
	{\it Let
	\[\mathbb{P}=\big\{t^{-k/p}f(t^{1/p})\mid f\in \IR[[t]], k\in\IN_0\mbox{ and }p\in\IN\big\}\]
	be the field of Puiseux series over $\IR$.
	\begin{itemize}
		\item[(1)]
		The volume algebra of $\mathbb{P}$ is the polynomial ring $\mathbb{P}[X]$ over $\mathbb{P}$.
		\item[(2)]
		Let $n\in\IN$.
		The maps
		\[\lambda_{\mathbb{P},n}:\chi^1_{\mathbb{P},n}\to \mathbb{P}[X]_{< n}, A\mapsto\lambda_{\mathbb{P},n}(A),\]
		and
		\[\mathrm{Int}_{\mathbb{P},n}:\ma{L}^1_{\mathbb{P},n}\to \mathbb{P}[X]_{\leq n}, f\mapsto \int_{\mathbb{P}^n}f\,d\lambda_{\mathbb{P},n},\]
		are surjective.
	\end{itemize}}
{\bf Proof:}

\vs{0.1cm}
	(1): The field $\IP$ of Puiseux series can be made canonically into a model of $T_\an$ (compare with Fact 1.5). So Theorem 3.18 gives (1).
	
	\vs{0.2cm}
	\noi (2): By Proposition 2.7 it is enough to deal with the integral. Let $n\in\IN$. We have that $\mathrm{Int}_{\IP,n}\big(\ma{L}^1_{\IP,n}\big)\subset \IP[X]_{\leq n}$ by Remark 3.16(2).
	That equality holds follows from Properties 2.3, Proposition 2.7 and the observation that
	\[\int_1^{t^{-1}}\cdots\int_1^{t^{-1}}\frac{d\lambda_{\IP,n}(x)}{x_1\cdot\ldots\cdot x_n}=\big(\log(t^{-1})\big)^n=X^n\]
	for all $n\in\IN$.
	\hfill$\Box$
}

\subsection{Functoriality of the construction}

\noi Let $R\subset S$ be an extension of real closed fields with archimedean value groups. We assume that $R$ contains the reals.
For a semialgebraic subset $A$ of $R^n$ or a semialgebraic function $f:R^n\to R$ let $A_S$ and $f_S:S^n\to S$ be the
canonical lifting of $A$ respectively $f$ to $S$ (via quantifier elimination, see for example [6, Chapter 5]).

\vs{0.5cm}
{\bf 3.20 Remark}

\vs{0.1cm}
	{\it Let $(s,\sigma,\tau)$ be a Lebesgue datum for $R$. Then there is a Lebesgue datum $(s^*,\sigma^*,\tau^*)$ for $S$ extending $(s,\sigma,\tau)$.}

\vs{0.5cm}
{\bf 3.21 Proposition}

\vs{0.1cm}
{\it 
	Let $s$ be a section for $R$ and let $s^*$ be a section for $S$ extending $s$.
	\begin{itemize}
		\item[(1)]
		Let $A\subset R^n$ be semialgebraic. Then $\lambda^{s^*}_{S,n}(A_S)=\lambda^s_{R,n}(A)$.
		\item[(2)]
		Let $f:R^n\to R$ be semialgebraic. Then $f$ is integrable over $R^n$ if and only if $f_S$ is integrable over $S^n$.
		If this holds then $\int f_S\,d\lambda^{s^*}_{S,n}=\int f\,d\lambda^s_{R,n}$.
	\end{itemize}}
{\bf Proof:}

\vs{0.1cm}
	This is evident from the construction of the measure and the integral.
\hfill$\Box$

\vs{0.5cm}
\noi In  particular we obtain that if the semialgebraic set or the semialgebraic function is defined over the reals then the value of the measure and the integral is the usual one obtained by measuring and integrating on the reals. We generalize this.

\vs{0.5cm}
\noi Let $R$ be a real closed field with archimedean value group that contains the reals.
The {\bf standard part map} $\mathrm{\bf st}:R\to \IR\cup\{\infty\}$ is defined as follows (see for example Van den Dries [20] or  Ma\v{r}\'{i}kov\'{a} [38]): Let $a\in R$. If $a$ is bounded (see Section 1.2) then $\mathrm{\bf st}(a)$ is the unique real number $x$ such that $x-a$ is infinitesimal. If $a$ is not bounded then $\mathrm{\bf st}(a)=\infty$. The standard part map for tuples is defined componentwise.

\noi We are interested in the behaviour of Lebesgue measure with respect to the standard map.
Similarly to above, the standard part map $\mathrm{\bf st}:\ma{R}[X]\to \IR\cup\{\infty\}$ is defined.

\noi We choose an arbitrary section for $R$.

\vs{0.5cm}
{\bf 3.22 Remark}	

\vs{0.1cm}				
	{\it Assume that $R\neq \IR$.
	For $n\geq 2$ there is an integrable semialgebraic subset $A$ of $R^n$ such that $\mathrm{\bf st}\big(\lambda_{R,n}(A)\big)\neq \lambda_n\big(\mathrm{\bf st}(A)\cap \IR^n\big)$.}

\vs{0.1cm}
{\bf Proof:}

\vs{0.1cm}
	We deal with the case $n=2$, the higher dimensional examples are constructed in a completely analogous way.
	Let $\varepsilon\in R_{>0}$ be infinitesimal.
	Let $A:=[-\varepsilon,\varepsilon]\times [-1/\varepsilon,1/\varepsilon]$. Then $\lambda_{R,2}(A)=4$ by Property 2.3(5) and so
	$\mathrm{\bf st}\big(\lambda_{R,2}(A)\big)=4$.
	We have $\mathrm{\bf st}(A)\cap\IR^2=\{0\}\times \IR$, hence $\lambda_2\big(\mathrm{\bf st}(A)\cap\IR^2\big)=0$.
\hfill$\Box$

\vs{0.5cm}
{\bf 3.23 Definition}

\vs{0.1cm}	
	A subset $A$ of $R^n$ is called $\IR$-bounded if there is some $C\in\IR_{>0}$ such that $|x|\leq C$ for all $x\in A$.

\vs{0.5cm}
\noi Note that a subset $A$ of $R^n$ is $\IR$-bounded if and only if $\mathrm{\bf st}(A)\subset \IR^n$.
Note also that, given $a\in R^n$, the set $\{a\}$ is $\IR$-bounded if and only if $a$ is bounded.

\vs{0.5cm}
{\bf 3.24 Remark}

\vs{0.1cm}
{\it 	Let $A\subset R^n$ be semialgebraic and $\IR$-bounded. Then $\mathrm{\bf st}(A)$ is a semialgebraic subset of $\IR^n$ and $\dim\big(\mathrm{\bf st}(A)\big)\leq \dim(A)$.}

\vs{0.1cm}
{\bf Proof:}

\vs{0.1cm}
	The first assertion follows from the Marker-Steinhorn theorem, see for example
	[20, Section 8]. The second one follows from [20, Proposition 9.3].
	\hfill$\Box$

\vs{0.5cm}
{\bf 3.25 Proposition}

\vs{0.1cm}
	{\it Let $A$ be a semialgebraic subset of $R^n$ that is $\IR$-bounded. Then \[\mathrm{\bf st}\big(\lambda_{R,n}(A)\big)=\lambda_n\big(\mathrm{\bf st}(A)\big).\]}

\vs{0.1cm}
{\bf Proof:}

\vs{0.1cm}
It holds that $\mathrm{\bf st}\big(A\triangle \mathrm{\bf st}(A)_R\big)\subset \mathrm{\bf st}(\partial A)$ where $\triangle$ denotes the symmetric difference, see [20, Lemma 10.1].
	
	\vs{0.2cm}
	\noi {\it Claim:} $\mathrm{\bf st}\big(\lambda_{R_\an,n}(A\triangle \mathrm{\bf st}(A)_R)\big)=0$.
	
	\vs{0.1cm}
	\noi {\it Proof of the claim:}
	We argue as in the proof of [20, Lemma 10.2]. Let $B:=\mathrm{\bf st}(\partial A)$. We have that $B$ is compact with $\dim(B)<n$ by Remark 3.24. Hence $\lambda_n(B)=0$.
	For $k\in\IN$ let $B_k:=\{x\in\IR^n\mid \mathrm{dist}(x,B)\leq 1/k\}$.
	We have that $B_k\searrow B$. Since $B_1$ is compact and hence $\lambda_n(B_1)<\infty$ we get by $\sigma$-continuity from above that $\lim_{k\to\infty}\lambda_n(B_k)=0$.
	Since $\mathrm{\bf st}\big(A\triangle \mathrm{\bf st}(A)_R\big)\subset B$ we get that $A\triangle \mathrm{\bf st}(A)_R\subset (B_k)_R$ for all
	$k\in\IN$.
	We obtain by Property 2.3(2) and Proposition 3.21 that
	\[\lambda_{R,n}\big(A\triangle (\mathrm{\bf st}(A))_R\big)\leq\lambda_{R,n}\big((B_k)_R\big)=\lambda_n(B_k).\]
	This shows that $\lambda_{R,n}\big(A\triangle \mathrm{\bf st}(A)_R\big)$ is infinitesimal.
	\hfill$\Box_{\mathrm{Claim}}$
	
	\vs{0.2cm}
	\noi Applying again Proposition 3.21, the claim gives that
	\begin{eqnarray*}
	\big\vert\lambda_{R,n}\big(A\big)-\lambda_n\big(\mathrm{\bf st}(A)\big)\big\vert&=&\big\vert\lambda_{R,n}\big(A\big)-\lambda_{R,n}\big(\mathrm{\bf st}(A)_R)\big\vert\\
	&\leq&\lambda_{R,n}\big(A\triangle \mathrm{\bf st}(A)_R\big)
	\end{eqnarray*}
	is infinitesimal. This shows the assertion.
\hfill$\Box$

\subsection{Extension to models of $T_\an$}

\noi Let $R$ be a real closed field with archimedean value group.
We assume that $R$ is a model of $T_\an$ (hence, it contains the reals). We write $R_\an$ when we want to stress that $R$ is considered as a model of $T_\an$ and not just as a real closed field.

\vs{0.5cm}
\noi Note that the results of Comte et al. (see Section 1.5) are formulated for globally subanalytic sets. Using Proposition 1.6, one can literally translate the constructions of the previous sections to obtain in the general case, given a section $s$ for $R$,
the {\bf analytic Lebesgue measure}
\[\lambda_{R_\an,n}=\lambda^s_{R_\an,n}:\big\{\mbox{globally subanalytic subsets of } R^n\big\}\to \ma{R}[X]\cup\{+\infty\}\]
and
the {\bf analytic Lebesgue integral}
\[\mathrm{Int}_{R_\an,n}:\ma{L}^1_{R_\an,n}\to \ma{R}[X]\]
{\bf on $R^n$
	with respect to the section}, and in the case that $\Gamma\cong \IQ$,
{\bf the analytic Lebesgue measure}
\[\lambda_{R_\an,n}:\big\{\mbox{globally subanalytic subsets of } R^n\big\}\to R[X]\cup\{+\infty\}\]
and {\bf the analytic Lebesgue integral}
\[\mathrm{Int}_{R_\an,n}:\ma{L}^1_{R_\an,n}\to R[X]\]
{\bf on $R^n$}
such that the above results hold, replacing semialgebraic by globally subanalytic.

\vs{0.5cm}
{\bf 3.26 Remark}

\vs{0.1cm}
{\it 
	The analytic Lebesgue measure and the analytic Lebesgue integral extend the semialgebraic Lebesgue measure and the semialgebraic Lebesgue integral (with respect to a section).}

\section{Constructible functions}

\noi The definition of constructible functions on $\IR_{\an,\exp}$ from the end of Section 1.5 can be naturally generalized to an arbitrary model of $T_{\an,\exp}$, in particular to the field $\ma{S}:=\IR((t))^{\mathrm{LE}}$ of LE-series.

\vs{0.5cm}
{\bf 4.1 Definition}

\vs{0.1cm}
A function $\ma{S}^n\to \ma{S}$ is called {\bf $\ma{S}$-constructible} if it is a finite sum of finite products of globally subanalytic functions and logarithms of positive globally subanalytic functions on $\ma{S}^n$.\\
\noi Similarly, we define an {\bf $\ma{S}$-constructible function on $A$} where $A$ is a globally subanalytic subset of some $\ma{S}^n$.	
	
\vs{0.5cm}
\noi 
Note that an $\ma{S}$-constructible function is definable in the $\ma{L}_{\an,\exp}$-structure $\ma{S}$.\\
By $\Omega_{\ma{S},n}$ we denote the ring of globally subanalytic functions $\ma{S}^n\to\ma{S}$. 
For $r>0$ we denote by $\ma{U}_\ma{S}(r)$ the set of all infinitely often differentiable globally subanalytic functions $f:]-r,r[\to \ma{S}$ such that $f(0)\neq 0$ and $|f(x)/f(0)-1|<1/2$ for all $x\in]-r,r[$. Let $\ma{U}_\ma{S}$ be the union of all $\ma{U}_\ma{S}(r)$.
Likewise for the $T_\an$-model $\IP$.\\
We have the following preparation result for unary constructible functions.

\vs{0.5cm}
{\bf 4.2 Proposition}

\vs{0.1cm}
	{\it Let $a,b\in \ma{S}$ with $a<b$ and let
	$f:]a,b[\to \ma{S}$ be $\ma{S}$-constructible. 
	Then there are $n\in \IN$ and $a=a_0<c_0<d_0<a_1<c_1<d_1<a_2<\ldots<a_{n-1}<c_{n-1}<d_{n-1}<a_n=b$ in $\ma{S}$ such that the following holds:
	\begin{itemize}
		\item[(1)]
		For every $j\in \{0,\ldots,n-1\}$ there
        is a finite subset $\ma{E}$ of $\IQ\times -\IN_0$, some $p\in\IN$ and for every $\sigma=(\sigma_1,\sigma_2)\in\ma{E}$ a function $u_\sigma\in\ma{U}_\ma{S}\big((d_j-a_j)^{1/p}\big)$ such that
        on $]a_j,d_j[$
        \[f(x)=\sum_{\sigma\in\ma{E}}u_\sigma(|x-a_j|^{1/p}) |x-a_j|^{\sigma_1}\big(\log |x-a_j|\big)^{-\sigma_2}.\]
        \item[(2)] 
        For every $j\in \{1,\ldots,n\}$ there
        is a finite subset $\ma{E}$ of $\IQ\times -\IN_0$, some $p\in\IN$ and for every $\sigma=(\sigma_1,\sigma_2)\in\ma{E}$ a function $u_\sigma\in\ma{U}_\ma{S}\big((a_j-c_{j-1})^{1/p}\big)$ such that
        on $]c_{j-1},a_j[$
        \[f(x)=\sum_{\sigma\in\ma{E}}u_\sigma(|x-a_j|^{1/p}) |x-a_j|^{\sigma_1}\big(\log |x-a_j|\big)^{-\sigma_2}.\]
     \end{itemize}}
{\bf Proof:}

\vs{0.1cm}
Let $g:]a,b[\to\ma{S}$ be globally subanalytic. Then it follows by the preparation theorem for polynomially bounded o-minimal structures (see van den Dries and Speissegger
[25, Lemma 2.2]) in connection with the characterization of globally subanalytic functions (see for example [37, p. 760]) that piecewise there are $\theta\in [a,b], \lambda\in\IQ, p\in\IN$ and $h\in\ma{U}_\ma{S}(r)$ for some $r>0$ such that $g(x)=|x-\theta|^\lambda h(|x-\theta|^{1/p})$ for all $x$ with either $0<x-\theta<r^p$ or $0<\theta-x<r^p$.
Let 
$c:=h(0)$. Assume that $c>0$. 
Applying the logarithm we get that $\log\big(g(x)\big)=\log\big(h(|x-\theta|^{1/p})\big)+\lambda \log |x-\theta|$.
By the definition of $\ma{U}_\ma{S}(r)$, we have that $\log\big(h(|x-\theta|^{1/p})\big)$ is defined and globally subanalytic on $]-r,r[$.
Taking the definition of a constructible function into account and applying once more the globally subanalytic preparation theorem to the functions of the latter type, we obtain the existence.
\hfill$\Box$

\vs{0.5cm}
\noi 
Note that the minus sign in the exponents of the logarithm is for technical reasons.

\vs{0.5cm}
{\bf 4.3 Corollary}

\vs{0.1cm}
{\it Let $f:\ma{S}_{>0}\to \ma{S}$ be $\ma{S}$-constructible. 
Then there is a finite subset $\ma{E}$ of $\IQ\times -\IN_0$, some $p\in\IN$ and for every $\sigma=(\sigma_1,\sigma_2)\in\ma{E}$ a function $u_\sigma\in\ma{U}_\ma{S}$ such that
\[f(x)=\sum_{\sigma\in\ma{E}}u_\sigma(x^{1/p}) x^{\sigma_1}\big(\log x\big)^{-\sigma_2}\]
for all sufficiently small positive $x\in\ma{S}$.}

\vs{0.5cm}
\noi We call the subset $\ma{E}$ of $\IQ\times -\IN_0$ from above {\bf a set of exponents}, $p$ {\bf a corresponding ramification index} and $(u_\sigma)_{\sigma\in\ma{E}}$ {\bf a corresponding coefficient tuple for $f$ at $0_+$}.
We equip $\IQ\times \IZ$ with the lexicographical ordering.

\vs{0.5cm}
{\bf 4.4 Remark}

\vs{0.1cm}
	{\it Let $(q,n)\in \IQ\times\IZ$. Then
	\[\lim_{x\to 0_+} x^q\big(\log x\big)^{-n}=
	\left\{\begin{array}{ccc}
	0,&& (q,n)>(0,0),\\
	1,&\mbox{iff}& (q,n)=(0,0),\\
	\infty,&&(q,n)<(0,0).
	\end{array}
	\right.\]}

\vs{0.5cm}
\noi We say that a function $f:\ma{S}_{>0}\to\ma{S}$ is ultimately zero at $0_+$ if $f(x)=0$ for all sufficiently small positive $x$.

\vs{0.5cm}
{\bf 4.5 Proposition}	

\vs{0.1cm}
	{\it Let $f:\ma{S}_{>0}\to \ma{S}$ be an $\ma{S}$-constructible function. Let $\ma{E}$ be a set of exponents, $p$ a corresponding ramification index and $(u_\sigma)_{\sigma\in \ma{E}}$ a corresponding coefficient tuple for $f$ at $0_+$.
	The following holds:
	\begin{itemize}
		\item[(1)]
		$f$ is ultimately zero at $0_+$ if and only if $\ma{E}=\emptyset$.
		\item[(2)] Assume that $f$ is not ultimately zero at $0_+$. Let $\mu_\ma{E}=(q_\ma{E},n_\ma{E}):=\min\ma{E}$.
		Then
		\[\lim_{x\to 0_+}\frac{f(x)}{x^{q_\ma{E}}\big(\log x\big)^{-n_\ma{E}}}=u_{\mu_\ma{E}}(0)\in\ma{S}\setminus\{0\}.\]
		In particular, $\mu_\ma{E}$ does not depend on the choice of $\ma{E}$, and we write $\mu_f=(q_f,n_f)$.
	\end{itemize}}
{\bf Proof:}

\vs{0.1cm}
	This follows from Remark 4.4.
\hfill$\Box$

\vs{0.5cm}
{\bf 4.6 Corollary}

\vs{0.1cm}
	{\it Let $f:\ma{S}_{>0}\to\ma{S}$ be an $\ma{S}$-constructible function that is not ultimately zero at $0_+$. 
	\begin{itemize}
		\item[(1)]
The function $f$ has a limit in $\ma{S}$ as $x$ tends to $0_+$ if and only if either $q_f>0$ or $q_f=n_f=0$.
\item[(2)] 
The function $f$ has the limit $0$ as $x$ tends to $0_+$ if and only if $q_f>0$.
\end{itemize}}

\vs{0.5cm}
\noi
We obtain the following quasianalyticity result for $\ma{S}$-constructible functions.

\vs{0.5cm}
{\bf 4.7 Proposition}

\vs{0.1cm}
{\it Let $f:\ma{S}^n\to \ma{S}$ be $\ma{S}$-constructible, and suppose that $f$ is $C^\infty$ in a neighbourhood of $a\in \ma{S}^n$ such that the Taylor series of $f$ at $a$ vanishes. Then $f=0$ on a neighbourhood of $a$.}

\vs{0.1cm}
{\bf Proof:}

\vs{0.1cm}
{\bf Special case:} We assume that $n=1$.
Applying a translation we may assume that $a=0$.
By symmetry it is enough to show that $f$ vanishes for all sufficiently small positive $x$.
Let $\ma{E}$ be a set of exponents for $f$ at $0_+$.
If $\ma{E}=\emptyset$ we are done by Proposition 4.5(1).
Assume that $\ma{E}\neq \emptyset$.
Let $\mu_f=(q_f,n_f)=\min\ma{E}$ be as in Proposition 4.5(2).
Since 
$\lim_{x\to 0}f(x)=0$ 
we see by Corollary 4.6(2) that $q_f>0$.
Choose $\lambda\in\IN$ with $\lambda>q_f$. By Proposition 4.5(2) we get that $f(x)>x^\lambda$ for all sufficiently small positive $x$.
But since the Taylor series of $f$ at $a$ vanishes we obtain by Taylor estimation that $\lim_{x\to 0}f(x)/x^n=0$ for all $n\in\IN$, contradiction.

\vs{0.2cm}
\noi{\bf General case:}
To prove the general case we strengthen the result of the special case:

\vs{0.2cm}
\noi{\bf Claim:}
Let $I$ be an open interval in $\ma{S}$ and let $g:I\to \ma{S}$ be an $\ma{S}$-constructible function that is $C^\infty$.
Assume that there exists $a\in I$ such that the Taylor series of $g$ at $a$ vanishes. Then $g=0$ on $I$.

\vs{0.1cm}
\noi{\bf Proof of the claim:}
Let
$B$ be the set of all $b\in I$ with $a<b$ such that  $g$ is identically zero on $[a,b[$. 
By the special case $B\neq\emptyset$. By o-minimality, $\sup B$ exists in $\ma{S}\cup\{\infty\}$. 
By the special case we obtain that $\sup B=\sup I$. Hence $g$ is identically zero on $I\cap [a,\infty[$. 
By symmetry we see that $g$ is identically zero on $I$. 
\hfill$\Box_{\mathrm{Claim}}$

\vs{0.2cm}
\noi Let $r\in\ma{S}_{>0}$ such that $f$ is $C^\infty$ on $U:=\big\{x\in\ma{S}^n\mid |x-a|<r\big\}$.
For each $b\in\ma{S}^n$ with $|b|=r$ define $g_b:]-1,1[\to \ma{S}, s\mapsto f(a+bs)$.
Then $g_b$ is an $\ma{S}$-constructible function that is $C^\infty$. The Taylor series of $g_b$ vanishes in $0$. By the claim we get that $g_b=0$ for every $b$.
Hence $f=0$ on $U$.
\hfill$\Box$	

\vs{0.5cm}
\noi
We define constructible functions on the field $\IP$ of Puiseux series. These functions take values in the volume algebra $\IP[X]$.
Throughout the section we use that $\IP$ carries a canonical $\ma{L}_\an$-structure which makes it a model of $T_\an$ (compare with Fact 1.5).

\vs{0.5cm}
{\bf 4.8 Remark}

\vs{0.1cm}
{\it 	
We have that $\log x\in-\mathrm{ord}(x) X+\ma{O}_\IP\subset\IP[X]$ for all $x\in\IP_{>0}$.}

\vs{0.1cm}
{\bf Proof:}

\vs{0.1cm}
See Lemma 3.2(2).
\hfill$\Box$

\vs{0.5cm}
{\bf 4.9 Definition}

\vs{0.1cm}
A function $\IP^n\to \IP[X]$ is called {\bf $\IP$-constructible} if it is a finite sum of finite products of globally subanalytic functions and logarithms of positive globally subanalytic functions on $\IP^n$.\\
	\noi Similarly, we define a {\bf $\IP$-constructible function on $A$} where $A$ is a globally subanalytic subset of some $\IP^n$.

\vs{0.5cm}
\noi The relevance of constructible functions for integrating is given by the following.

\vs{0.5cm}
{\bf 4.10 Proposition}

\vs{0.1cm}
	{\it Let $f:\IP^{q+n}\to \IP$ be globally subanalytic. Then the following holds:
	\begin{itemize}
		\item[(1)]
		The set
		\[\mathrm{Fin}(f):=\Big\{t\in \IP^q\;\big\vert\, f_t\mbox{ is integrable}\Big\}\]
		is globally subanalytic.
		\item[(2)]
		There is a $\IP$-constructible function $h:\IP^q\to\IP[X]$ such that
		\[\int_{\IP^n}f_t\,d\lambda_{\IP,n}=h(t)\]
		for all $t\in \mathrm{Fin}(f)$.
	\end{itemize}}
{\bf Proof:}

\vs{0.1cm}
	This follows by doing Construction 2.6 with parameters.
\hfill$\Box$

\vs{0.5cm}
{\bf 4.11 Proposition}

\vs{0.1cm}
	{\it Let $a,b\in \IP$ with $a<b$ and let
	$f:\IP\to \IP[X]$ be $\IP$-constructible. 
	Then there are $n\in \IN$ and $a=a_0<c_0<d_0<a_1<c_1<d_1<a_2<\ldots<a_{n-1}<c_{n-1}<d_{n-1}<a_n=b$ in $\IP$ such that the following holds:
	\begin{itemize}
		\item[(1)]
		For every $j\in \{0,\ldots,n-1\}$ there
		is a finite subset $\ma{F}$ of $\IQ\times -\IN_0\times\IN_0$, some $p\in\IN$ and for every $\tau=(\tau_1,\tau_2,\tau_3)\in\ma{F}$ a function $v_\tau\in\ma{U}_\IP\big((d_j-a_j)^{1/p}\big)$ such that
		on $]a_j,d_j[$
		\[f(x)=\sum_{\tau\in\ma{F}}v_\tau(|x-a_j|^{1/p}) |x-a_j|^{\tau_1}\big(\log |x-a_j|\big)^{-\tau_2}X^{\tau_3}.\]
		\item[(2)] 
		For every $j\in \{1,\ldots,n\}$ there
		is a finite subset $\ma{F}$ of $\IQ\times -\IN_0\times\IN_0$, some $p\in\IN$ and for every $\tau=(\tau_1,\tau_2,\tau_3)\in\ma{F}$ a function $v_\sigma\in\ma{U}_\IP\big((a_j-c_{j-1})^{1/p}\big)$ such that
		on $]c_{j-1},a_j[$
		\[f(x)=\sum_{\tau\in\ma{F}}v_\tau(|x-a_j|^{1/p}) |x-a_j|^{\tau_1}\big(\log |x-a_j|\big)^{-\tau_2}X^{\tau_3}.\]
	\end{itemize}}
{\bf Proof:}

\vs{0.1cm}
Let $g:]a,b[\to\IP$ be globally subanalytic. Then it follows by the preparation theorem for polynomially bounded o-minimal structures (see van den Dries and Speissegger
[25, Lemma 2.2]) in connection with the characterization of globally subanalytic functions (see for example [37, p. 760]) that piecewise there are $\theta\in [a,b], \lambda\in\IQ, p\in\IN$ and $h\in\ma{U}_\IP(r)$ for some $r>0$ such that $g(x)=|x-\theta|^\lambda h(|x-\theta|^{1/p})$ for all $x$ with either $0<x-\theta<r^p$ or $0<\theta-x<r^p$.
Let 
$c:=h(0)$ and $\widetilde{h}:=h/c$. Assume that $c>0$. 
Applying the logarithm we get that $\log\big(g(x)\big)=\log(c)+\log\big(\widetilde{h}(|x-\theta|^{1/p})\big)+\lambda \log |x-\theta|$.
We have that $\log(c)\in \IR X+\ma{O}_\IP$ by Remark 4.8 and, by the definition of $\ma{U}_\IP(r)$, that $\log\big(\widetilde{h}(|x-\theta|^{1/p})\big)$ is defined and globally subanalytic on $]-r,r[$.
Taking the definition of a constructible function into account and applying once more the globally subanalytic preparation theorem to the functions of the latter type, we obtain the existence.
\hfill$\Box$

\vs{0.5cm}
{\bf 4.12 Corollary}

\vs{0.1cm}
	{\it Let $f:\IP_{>0}\to \IP[X]$ be $\IP$-constructible. 
	Then there is a finite subset $\ma{F}$ of $\IQ\times -\IN_0\times \IN_0$, some $p\in\IN$ and for every $\tau=(\tau_1,\tau_2,\tau_3)\in\ma{F}$ a function $v_\tau\in\ma{U}_\IP$ such that
	\[f(x)=\sum_{\tau\in\ma{F}}v_\tau(x^{1/p}) x^{\tau_1}\big(\log x\big)^{-\tau_2}X^{\tau_3}\]
	for all sufficiently small positive $x\in\IP$.}

\vs{0.5cm}
\noi We call the subset $\ma{F}$ of $\IQ\times -\IN_0\times\IN_0$ from above {\bf a set of exponents}, $p$ {\bf a corresponding ramification index} and $(v_\tau)_{\tau\in\ma{F}}$ {\bf a corresponding coefficient tuple for $f$ at $0_+$}.\\
\noi From Remark 4.8 one sees that $\log x$ and all of its positive powers do not have a limit as $x$ tends to $0_+$ in $\IP$.

\vs{0.2cm}
\noi An easy calculation gives the following:

\vs{0.5cm}
{\bf 4.13 Remark}

\vs{0.1cm}
	{\it Let $(q,n)\in \IQ\times\IZ$. Then, in $\IP(X)$, we have that
	\[\lim_{x\to 0_+} x^q\big(\log x\big)^{-n}=
	\left\{\begin{array}{ccc}
	0,&& q>0,\\
	1,&\mbox{iff}& q=n=0,\\
	\infty,&&q<0.
	\end{array}
	\right.\]}

\vs{0.5cm}
\noi Let $\pi:\IQ\times-\IN_0\times\IN_0\to \IQ\times-\IN_0$ be the projection onto the first two components. As above, we equip $\IQ\times \IZ$ with the lexicographical ordering.

\vs{0.5cm}
{\bf 4.14 Proposition}	

\vs{0.1cm}
	{\it Let $f:\IP_{>0}\to \IP[X]$ be a $\IP$-constructible function. Let $\ma{F}$ be a set of exponents, $p$ a corresponding ramification index and $(v_\tau)_{\tau\in\ma{F}}$ a corresponding coefficient tuple for $f$ at $0_+$.
	The following holds:
	  \begin{itemize}
		\item[(1)]
		$f$ is ultimately zero at $0_+$ if and only if $\ma{F}=\emptyset$.
		\item[(2)] Assume that $f$ is not ultimately zero at $0_+$. Let $\nu_\ma{F}=(q_\ma{F},n_\ma{F}):=\min\pi(\ma{F})$.
		Then
		\[\lim_{x\to 0_+}\frac{f(x)}{x^{q_\ma{F}}\big(\log x\big)^{-n_\ma{F}}}=\sum_{\pi(\tau)=\nu_\ma{F}}v_\tau(0)X^{\tau_3}\in \IP[X]\setminus\{0\}.\]		
		In particular, $\nu_\ma{F}$ does not depend on the choice of $\ma{F}$ and we write $\nu_f=(q_f,n_f)$.
	\end{itemize}}
{\bf Proof:}

\vs{0.1cm}
	This follows from Remark 4.13 and the transcendence of $X$ over $\IP$.
\hfill$\Box$

\vs{0.5cm}
{\bf 4.15 Corollary}

\vs{0.1cm}
{\it Let $f:\IP_{>0}\to \IP[X]$ be a $\IP$-constructible function that is not ultimately zero at $0_+$.
\begin{itemize}
	\item[(1)]
The function $f$ has a limit in $\IP[X]$ as $x$ tends to $0_+$ if and only if
either $q_f>0$ or $q_f=n_f=0$.
\item[(2)] The function $f$ has the limit $0$ as $x$ tends to $0_+$ if and only if
$q_f>0$. 
\end{itemize}}

\vs{0.2cm}
\noi We want to lift $\IP$-constructible functions to $\ma{S}$-constructible functions. Certainly, this cannot be done in a unique way.
But we want to define a canonical lifting. The idea is to use the definition of constructible functions.

\noi Let $g:\IP^n\to \IP$ be globally subanalytic.
By $g_\ma{S}:\ma{S}^n\to \ma{S}$ we denote the canonical lifting of $g$ from $\IP$ to $\ma{S}$ (as models of the theory $T_\an=\mathrm{Th}(\IR_\an)$).

\vs{0.5cm}
{\bf 4.16 Definition}

\vs{0.1cm}
	Let $f:\IP^n\to \IP[X]$ be $\IP$-constructible.
	\begin{itemize}
		\item[(a)]
		There are $r\in\IN$, a polynomial $Q\in \Omega_{\IP,n}[T_1,\ldots,T_r]$ and positive $\varphi_1,\ldots,\varphi_r\in\Omega_{\IP,n}$ such that $f(x)=Q\big(\log(\varphi_1(x)),\ldots,\log(\varphi_r(x))\big)$ for all $x\in \IP$.
		We call $\Delta:=(r,Q,\varphi)$ where $\varphi=(\varphi_1,\ldots,\varphi_r)$ a {\bf representation tuple for $f$}.
		\item[(b)]
		Let $\Delta=(r,Q,\varphi)$ be a representation tuple for $f$.
		We write
		$f_{\ma{S},\Delta}$ for the $\ma{S}$-constructible function $$Q_\ma{S}\big(\log((\varphi_1)_\ma{S}),\ldots,\log((\varphi_r)_\ma{S})\big)$$ on $\ma{S}^n$.
	\end{itemize}

\vs{0.2cm}
\noi Clearly $f_{\ma{S},\Delta}$ lifts $f$ to $\ma{S}$.

\vs{0.5cm}
{\bf 4.17 Proposition}

\vs{0.1cm}
{\it Let $f:\IP^n\to\IP[X]$ be a $\IP$-constructible function and let $\Delta$ be a representation tuple for $f$. If $f$ is identically zero on $\IP^n$ then $f_{\ma{S},\Delta}$ is identically zero on $\ma{S}^n$.}

\vs{0.1cm}
{\bf Proof:}

\vs{0.1cm}
Let $\Delta=(r,Q,\varphi)$. We can find a partition $\ma{C}$ of $\IP^n$ into globally subanalytic $C^\infty$-cells such that for each $C\in\ma{C}$, the restrictions to $C$ of $\varphi$ and of the globally subanalytic coefficient functions of the polynomial $Q\in \Omega_{\IP,n}[T_1,\ldots,T_r]$ are all $C^\infty$.
Fix $C\in\ma{C}$ and let $C_\ma{S}$ be the canonical lifting of $C$ to the $T_\an$-model $\ma{S}$. Note that $C=C_\ma{S}\cap \IP^n$. Set $Z:=\big\{x\in C_\ma{S}\mid f_{\ma{S},\Delta}(x)=0\big\}$. It suffices to prove that $Z=C_\ma{S}$.
The liftings of $\varphi$ and of all the globally subanalytic coefficient functions of the polynomial $Q$ to $\ma{S}$ are $C^\infty$ on $C_\ma{S}$.
So $f_{\ma{S},\Delta}$ is $C^\infty$ on $C_\ma{S}$. In particular, $f_{\ma{S},\Delta}$ is continuous on $C_\ma{S}$, so $Z$ is closed in $C_\ma{S}$.

\vs{0.2cm}
\noi{\bf Case 1:} $\dim C=n$.

\vs{0.1cm}
\noi{\bf Proof of Case 1:}
Then $C$ is open in $\IP^n$ and $C_\ma{S}$ is open in $\ma{S}^n$. Working in the model $\ma{S}$ of $T_{\an,\exp}$, we find a $C^\infty$-cell decomposition $\ma{D}$ of $C_\ma{S}$ (in the sense of $\ma{L}_{\an,\exp}$) which is compatible with $Z$ and stratifies $C_\ma{S}$.
Let 
$$\ma{D}_1:=\big\{D\in\ma{D}\mid D\subset Z\big\},\ma{D}_2:=\big\{D\in \ma{D}\mid D\cap \IP^n\neq\emptyset\big\},$$
$$\ma{D}_3:=\big\{D\in\ma{D}\mid D\mbox{ open}\big\}.$$ 
Since $f_{\ma{S},\Delta}(x)=0$ for all $x\in C=C_\ma{S}\cap\IP^n$ it follows by the compability of $\ma{D}$ with $Z$ that $\ma{D}_2\subset \ma{D}_1$.

\vs{0.2cm}
\noi {\bf Claim:}
There exists an open cell $D\in\ma{D}$ such that $D\cap \IP^n\neq\emptyset$; i.e. $\ma{D}_2\cap\ma{D}_3\neq\emptyset$.

\vs{0.1cm}
\noi{\bf Proof of the claim:}
We prove the claim by induction on $n$. First assume that $n=1$. Then $\ma{D}$ is a finite set of points and intervals. Since $C=C_\ma{S}\cap\IP$ is an infinite subset of $C_\ma{S}$, there must exist some open interval $D\in\ma{D}$ such that $D\cap \IP\neq\emptyset$.\\
Next we assume that $n>1$ and that the claim holds with $n-1$ in place of $n$. By applying the base case of the induction to the open interval $\pi_1(C_\ma{S})=\big(\pi_1(C)\big)_\ma{S}$ (where $\pi_1(x_1,\ldots,x_n)=x_1$), we may choose an open interval $B\in \{\pi_1(D)\mid D\in\ma{D}\}$ such that $B\cap\IP\neq \emptyset$. Let $b\in B\cap \IP$.
Since the fiber $(C_\ma{S})_b=(C_b)_\ma{S}$ is an open cell in $\ma{S}^{n-1}$ and the set of fibers $\big\{D_b\mid D\in \ma{D}\mbox{ s.t. }\pi_1(D)=B\big\}$ is a cell decomposition of $(C_\ma{S})_b$, the induction hypothesis supplies a $D\in\ma{D}$ such that $\pi_1(D)=B$, $D_b$ is open in $\ma{S}^{n-1}$, and $D_b\cap \IP^{n-1}\neq \emptyset$. It follows that $D$ is an open cell and $D\cap \IP^n\neq\emptyset$.
\hfill$\Box_{\mathrm{Claim}}$

\vs{0.2cm}
\noi
Set $\Sigma:=\ma{D}_1\cap\ma{D}_3$. We have that $\Sigma\supset \ma{D}_2\cap \ma{D}_3$ where the latter is nonempty by the claim. Let $Y$ be the closure of the union of the cells that are elements of $\Sigma$.
Since $Z$ is closed we get that $Y\subset Z$. So to prove that $Z=C_\ma{S}$, it suffices to prove that $Y=C_\ma{S}$.
Note that $Y$ is a union of cells in $\ma{D}$ since $\ma{D}$ stratifies $C_\ma{S}$.
Let $a\in Y$. Fix $D\in\Sigma$ such that $a$ is in the closure of $D$. Since $f_{\ma{S},\Delta}=0$ on the open set $D$ and $f_{\ma{S},\Delta}$ is $C^\infty$ on $C_\ma{S}$, the Taylor series of $f_{\ma{S},\Delta}$ at $a$ vanishes.
By Proposition 4.7 we obtain that $f_{\ma{S},\Delta}=0$ on a neighbourhood $U$ of $a$ in $C_\ma{S}$.
It follows that $a$ is in the interior of $Y$, and this shows that $Y$ is open.
So in summary,  $Y$ is a nonempty $\ma{L}_{\an,\exp}$-definable set that is both open and closed in $C_\ma{S}$. Since $C_\ma{S}$ is definably connected, in the sense of $\ma{L}_{\an,\exp}$, it follows that 
$Y=C_\ma{S}$.

\vs{0.2cm}
\noi{\bf Case 2:} $\dim C<n$.

\vs{0.1cm}
\noi{\bf Proof of Case 2:}
The cell $C$ is the graph of a $C^\infty$ globally subanalytic section of the projection of $C$ onto an open cell in a lower dimensional space, so by pulling back $f$ and its representation tuple $\Delta$ by that section we arrive at Case 1.
\hfill$\Box$

\vs{0.5cm}
{\bf 4.18 Theorem}

\vs{0.1cm}
	{\it Let $f:\IP^n\to \IP[X]$ be $\IP$-constructible. Let $\Delta,\widehat{\Delta}$ be representation tuples for $f$.
	Then $f_{\ma{S},\Delta}=f_{\ma{S},\widehat{\Delta}}$.}

\vs{0.1cm}
{\bf Proof:}

\vs{0.1cm}
Let $\Delta=(Q,r,\varphi),\widehat{\Delta}=(\widehat{Q},\widehat{r},\widehat{\varphi})$ be representation tuples for $f$. Since $f-f=0$ we see that $\Delta^*:=(Q-\widehat{Q},(\varphi,\widehat{\varphi}),r+\widehat{r})$ is a representation tuple for the zero function on $\IP^n$ whose lifting to $\ma{S}$ with respect to $\Delta^*$ is given by $f_{\ma{S},\Delta}-f_{\ma{S},\widehat{\Delta}}$. By Proposition 4.17 we obtain that $f_{\ma{S},\Delta}-f_{\ma{S},\widehat{\Delta}}$ is identically zero  on $\ma{S}^n$, so $f_{\ma{S},\Delta}=f_{\ma{S},\widehat{\Delta}}$.
\hfill$\Box$

\vs{0.5cm}
\noi Consequently, we write $f_\ma{S}$ for $f_{\ma{S},\Delta}$ where $\Delta$ is a presentation tuple for $f$ and call it the {\bf canonical lifting of $f$}. This, of course, generalizes the case that $f$ is globally subanalytic.

\noi Note that in the case $n=1$, by o-minimality, $\lim_{x\to 0_+,x\in\ma{S}}f_{\ma{S}}(x)$ exists in $\ma{S}\cup\{\pm\infty\}$.

\vs{0.5cm}
{\bf 4.19 Proposition}

\vs{0.1cm}
	{\it Let $f:\IP_{>0}\to\IP[X]$ be $\IP$-constructible. 
	\begin{itemize}
		\item[(1)] $f$ is ultimately zero at $0_+$ if and only if $f_\ma{S}$ is ultimately zero at $0_+$.
    \item[(2)] 
    Assume that $f$ is not ultimately zero at $0_+$.
	Let $\ma{F}$ be a set of exponents, $p$ a corresponding ramification index and $(v_\tau)_{\tau\in \ma{F}}$ be a corresponding coefficient tuple for $f$ at $0_+$.
	For $\sigma\in \pi(\ma{F})$ set $u_\sigma:=\sum_{\pi(\tau)=\sigma}(v_\tau)_\ma{S} X^{\tau_3}$.
	Then $\pi(\ma{F})$ is a set of exponents, $p$ a corresponding ramification index and  $(u_\sigma)_{\sigma\in\pi(\ma{F})}$ a corresponding coefficient tuple for $f_\ma{S}$ at $0_+$. In particular, $\nu_f=\mu_{f_\ma{S}}$.
	\end{itemize}}
{\bf Proof:}

\vs{0.1cm}
	(1): This follows immediately from Proposition 4.17.
	
	\vs{0.2cm}
	\noi (2): From
	\[f(x)=\sum_{\tau\in\ma{F}}v_\tau(x^{1/p}) x^{\tau_1}\big(\log x\big)^{\tau_2}X^{\tau_3}\]
	for all sufficiently small positive $x\in \IP$ we obtain
	\[f_\ma{S}(x)=\sum_{\sigma\in\pi(\ma{F})}u_\sigma(x^{1/p})x^{\sigma_1}\big(\log x\big)^{\sigma_2}\]
	for all sufficiently small positive $x\in\ma{S}$.
	By the transcendence of $X$ over $\IP$ we see that $u_\sigma\in \ma{U}_\ma{S}$ for every $\sigma\in \pi(\ma{F})$. This gives the proposition.
\hfill$\Box$

\vs{0.5cm}
{\bf 4.20 Theorem}	

\vs{0.1cm}
{\it 	
	Let $f:\IP\to \IP[X]$ be $\IP$-constructible. The following are equivalent:
	\begin{itemize}
		\item[(i)] $\lim_{x\to 0_+,x\in \IP}f(x)$ exists in $\IP[X]$,
		\item[(ii)] $\lim_{x\to 0_+,x\in\ma{S}}f_{\ma{S}}(x)$ exists in $\ma{S}$.
	\end{itemize}
	If this holds then
	$\lim_{x\to 0_+,x\in\ma{S}}f_{\ma{S}}(x)=\lim_{x\to 0_+,x\in \IP}f(x)$.}

\vs{0.1cm}
{\bf Proof:}

\vs{0.1cm}
	By Proposition 4.19(1) we can assume that $f$ is not ultimately zero at $0_+$.
	By Proposition 4.19(2) we have that 
	$\nu_f=\mu_{f_\ma{S}}$.
	By Corollary 4.6(1) and Corollary 4.15(1) we obtain that (i) and (ii) are equivalent; applying Proposition 4.19(2) in connection with Proposition 4.5(2) and Proposition 4.14(2) we get equality of the limits in the case of existence.
	\hfill$\Box$

\vs{0.5cm}
\noi A function $f:\IP\to \IP[X]$ is, as usual, called differentiable at $x\in \IP$ if the limit
\[f'(x):=\lim_{y\to x}\frac{f(y)-f(x)}{y-x}\]
exists in $\IP[X]$.
Partial derivatives of functions $\IP^n\to \IP[X]$ are defined arccordingly.

\vs{0.5cm}
{\bf 4.21 Theorem}

\vs{0.1cm}
	{\it Let $f:\IP\to\IP[X]$ be $\IP$-constructible and let $x\in \IP$.
	The following are equivalent:
	\begin{itemize}
		\item[(i)] $f$ is differentiable at $x$.
		\item[(ii)] $f_\ma{S}$ is differentiable at $x$.
	\end{itemize}
	If this holds then $f'(x)=(f_\ma{S})'(x)$.}

\vs{0.1cm}
{\bf Proof:}

\vs{0.1cm}
	Apply Theorem 4.20 to the functions
	\[\IP_{>0}\to \IP[X], t\mapsto \frac{1}{t}\big(f(x+t)-f(x)\big)\]
	and
	\[\IP_{>0}\to \IP[X], t\mapsto \frac{1}{-t}\big(f(x-t)-f(x)\big).\]
\hfill$\Box$

\vs{0.5cm}
{\bf 4.22 Example}

\vs{0.1cm}
	{\it The logarithm $\log:\IP_{>0}\to \IP[X]$ is continuously differentiable with $(\log x)'=1/x$ for all $x\in \IP_{>0}$.}

\vs{0.5cm}
\noi Let $I$ be an open interval of $\IP$.

\vs{0.5cm}
{\bf 4.23 Corollary}

\vs{0.1cm}
{\it Let $f:I\to \IP[X]$ be $\IP$-constructible.
The following are equivalent:
\begin{itemize}
	\item[(i)] $f:I\to \IP[X]$ is differentiable.
	\item[(ii)] $f_\ma{S}:I_\ma{S}\to \ma{S}$ is differentiable. 
\end{itemize}}
{\bf Proof:}

\vs{0.1cm}
Let $(r,Q,\varphi)$ be a representation tuple for $f$. By o-minimality, $\varphi$ and the globally subanalytic coefficient functions of $Q$ are differentiable outside a finite subset $A$ of $I$. Their liftings to $\ma{S}$ are differentiable on $I_\ma{S}\setminus A$.
Since the function $\log:\ma{S}_{>0}\to \ma{S}$ is differentiable we get, by the usual chain rule and the product rule, that $f_\ma{S}$ is differentiable on $I_\ma{S}\setminus A$. We obtain the assertion of the corollary by Theorem 4.21.
\hfill$\Box$

\vs{0.5cm}
{\bf 4.24 Theorem}

\vs{0.1cm}
{\it Let $f:I\to \IP[X]$ be $\IP$-constructible and continuously differentiable. 
Then $f'$ is $\IP$-constructible and $(f')_\ma{S}=(f_\ma{S})'$.}

\vs{0.1cm}
{\bf Proof:}

\vs{0.1cm}
By Corollary 4.23 we know that $f_\ma{S}$ is differentiable.  
Let $(r,Q,\varphi)$ be a representation tuple for $f$. By o-minimality, $\varphi$ and the globally subanalytic coefficient functions of $Q$ are differentiable with globally subanalytic derivatives outside a finite subset $A$ of $I$.
Their liftings to $\ma{S}$ are differentiable on $I_\ma{S}\setminus A$. 
Let $g:=f|_{I\setminus A}$.
Applying Example 4.22 resp. the fact that $\log:\ma{S}_{>0}\to \ma{S}$ is continuously differentiable we get, by the usual chain rule and the product rule, that $g'$ is $\IP$-constructible on $I\setminus A$, that $(g_\ma{S})'$ is 
$\ma{S}$-constructible on $I_\ma{S}\setminus A$ and that $(g')_\ma{S}(x)=(g_\ma{S})'(x)$ for all $x\in I_\ma{S}\setminus A$.
Since $A$ is finite we get that $f'$ is $\IP$-constructible. 
For $a\in A\subset \IP$ we obtain from the condition that $f'$ is continuous and by Theorem 4.20 that
\begin{eqnarray*}
(f')_\ma{S}(a)&=&f'(a)=\lim_{x\to a, x\in\IP}f'(x)=\lim_{x\to a,x\in\IP} g'(x)\\
&=&\lim_{x\to a,x\in\ma{S}}(g')_\ma{S}(x)=\lim_{x\to a,x\in\ma{S}} (g_\ma{S})'(x)\\
&=&\lim_{x\to a,x\in\ma{S}} (f_\ma{S})'(x).
\end{eqnarray*}
The function $f_\ma{S}$ is definable in the model $\ma{S}$ of the o-minimal theory $T_{\an,\exp}$, and differentiable. Hence it is continuously differentiable. Therefore $\lim_{x\to a,x\in\ma{S}}(f_\ma{S})'(x)=(f_\ma{S})'(a)$ and we are done. 
\hfill$\Box$

\vs{0.5cm}
{\bf 4.25 Definition}

\vs{0.1cm}
	Let $f:I\to \IP$ be globally subanalytic.
	A $\IP$-constructible function $F:I\to \IP[X]$ is called an {\bf antiderivative of $f$}
	if $F$ is differentiable with $F'=f$.

\vs{0.5cm}
{\bf 4.26 Proposition}	

\vs{0.1cm}	
	{\it Let $f:I\to \IP$ be a globally subanalytic function that is continuous.
	Let $F,\widehat{F}:I\to \IP[X]$ be antiderivatives of $f$. Then there is some $c\in \IP[X]$ such that $\widehat{F}=F+c$.}

\vs{0.1cm}
{\bf Proof:}

\vs{0.1cm}
	Let $G:=F_{\ma{S}}$ and let $\widehat{G}:=\widehat{F}_{\ma{S}}$.
	By Corollary 4.23, $G$ and $\widehat{G}$ are differentiable. By Theorem 4.24 we have that
	\[G'=(F')_\ma{S}=f_\ma{S}=(\widehat{F}')_\ma{S}=\widehat{G}'.\]
    Since $G$ and $\widehat{G}$ are definable in the model $\ma{S}$ of the o-minimal theory $T_{\an,\exp}$ we can apply the mean value property and get that there is a constant $c\in\ma{S}$ such that $G=\widehat{G}+c$. Hence $\widehat{F}=F+c$ and we see that
	$c\in \IP[X]$.
\hfill$\Box$

\vs{0.5cm}
\noi We end the section with an application of the preparation result in Proposition 4.11 which we need in Section 6.

\vs{0.5cm}
{\bf 4.27 Theorem}

\vs{0.1cm}
	{\it Let $f:I\to \IP[X]$ be a $\IP$-constructible function that is infinitely often differentiable.
	Let $K$ be a closed and bounded subinterval of $I$. Then there are $N\in\IN_0$ and uniquely determined globally subanalytic functions $h_0,\ldots,h_N$ that are infinitely often differentiable on an open neighbourhood of $K$ such that $f|_K=\sum_{j=0}^Nh_jX^j$.}

\vs{0.1cm}
{\bf Proof:}

\vs{0.1cm}
	Let $a,b\in I$ with $a<b$ and $K\subset ]a,b[$. 
	By Proposition 4.11 
	there are $n\in \IN$ and $a=a_0<c_0<d_0<a_1<c_1<d_1<a_2<\ldots<a_{n-1}<c_{n-1}<d_{n-1}<a_n=b$ in $\IP$ such that the following holds:
	\begin{itemize}
		\item[(1)]
		For every $j\in \{0,\ldots,n-1\}$ there
		is a finite subset $\ma{F}$ of $\IQ\times -\IN_0\times\IN_0$, some $p\in\IN$ and for every $\tau=(\tau_1,\tau_2,\tau_3)\in\ma{F}$ a function $v_\tau\in\ma{U}_\IP\big((d_j-a_j)^{1/p}\big)$ such that
		on $]a_j,d_j[$
		\[f(x)=\sum_{\tau\in\ma{F}}v_\tau(|x-a_j|^{1/p}) |x-a_j|^{\tau_1}\big(\log |x-a_j|\big)^{-\tau_2}X^{\tau_3}.\]
		\item[(2)] 
		For every $j\in \{1,\ldots,n\}$ there
		is a finite subset $\ma{F}$ of $\IQ\times -\IN_0\times\IN_0$, some $p\in\IN$ and for every $\tau=(\tau_1,\tau_2,\tau_3)\in\ma{F}$ a function $v_\sigma\in\ma{U}_\IP\big((a_j-c_{j-1})^{1/p}\big)$ such that
		on $]c_{j-1},a_j[$
		\[f(x)=\sum_{\tau\in\ma{F}}v_\tau(|x-a_j|^{1/p}) |x-a_j|^{\tau_1}\big(\log |x-a_j|\big)^{-\tau_2}X^{\tau_3}.\]
	\end{itemize}
	Let $j\in \{0,\ldots,n-1\}$. We consider the interval $]a_j,d_j[$. Without restriction we may assume that $a_j=0$. 
	Let $g:=f(x^p)$. Then $g$ is $C^\infty$ on a neighbourhood of $0$ and
	\[g(x)=\sum_{\tau\in\ma{F}}p^{-\tau_2}v_\tau(x) x^{p\tau_1}\big(\log x\big)^{-\tau_2}X^{\tau_3}\]
	for all sufficiently small positive $x$.
	By Proposition 4.14(2) and the fact that $g$ is $C^\infty$ at $0$ we obtain that 
	$\nu_f=\min\pi(\ma{F})\in \frac{1}{p}\IN_0\times \{0\}$ if $\ma{F}\neq\emptyset$.
	Subtracting $\sum_{\pi(\tau)=\nu_f}p^{-\tau_2}v_\tau(x) x^{p\tau_1}X^{\tau_3}$ from $g$ and repeating this argument for the new function which has again a $C^\infty$-extension to $0$ we get that
	$\pi(\ma{F})\subset \frac{1}{p}\IN_0\times\{0\}$, i.e. no logarithmic terms occur.
	Repeating this for the intervals of the form $]c_{j-1},a_j[$ 
	we find $N\in\IN_0$ and globally subanalytic functions $h_0,\ldots,h_N:]a,b[\to \IP$ such that
	$f=\sum_{j=0}^Nh_jX^j$ on $]a,b[$. Since $X$ is transcendental over $\IP$ we get that $h_0,\ldots,h_N$ are $C^\infty$ on $]a,b[$ and that these functions are uniquely determined.
\hfill$\Box$

\section{Main theorems of integration}

\noi In this section we establish the transformation formula, Lebesgue's theorem on dominated convergence and the fundamental theorem of calculus for semialgebraic and analytic integration.

\noi As explained in the introduction, we deal with the field of Puiseux series. The results can be immediately translated to a model of $T_\an$ with archimedean value group, independently of the choice of the section.
For the general situation, i.e. for a real closed with archimedean value group that contains the reals, the results, again independently of the chosen section, can be also obtained by a more technical set up. We have refrained from this to concentrate on the key ideas.

\noi One can naturally realize the analytic, respectively, the semialgebraic, Lebesgue measure and integral on the field of Puiseux series by the Lebesgue datum $(s,\sigma,\tau)$, where $s$ maps $q\in\IQ$ to $t^q$ and $\sigma:\IP\hookrightarrow  \ma{P}=\IR((t^\IQ))$ and $\tau: \IQ\hookrightarrow\IR$ are the inclusions, respectively.

\noi We formulate the results below for the globally subanalytic setting. The corresponding statements for the semialgebraic one are then automatically included.

\subsection{The transformation formula}

\noi
{\bf 5.1 Theorem} (Transformation formula)
	
\vs{0.1cm}
{\it Let $U,V\subset \IP^n$ be globally subanalytic sets that are open and let $\varphi:U\to V$ be a globally subanalytic $C^1$-diffeomorphism.
	Let $f:V\to R$ be globally subanalytic.
	Then
	$f$ is integrable over $V$ if and only if $(f\circ\varphi)\big\vert \det(D_\varphi)\big\vert$ is integrable over $U$, and in this case
	\[\int_V f\,d\lambda_{\IP,n}=\int_U(f\circ\varphi)\big\vert\det(D_\varphi)\big\vert\,d\lambda_{\IP,n}.\]}

\vs{0.1cm}
{\bf Proof:}
	
	\vs{0.1cm}
	The usual transfer argument of Section 2 does the job.
\hfill$\Box$

\subsection{Lebesgue's theorem on dominated convergence and the fundamental theorem of calculus}

\noi Lebesgue's theorem on dominated convergence and the fundamental theorem of calculus involve limits with respect to the raw data. By simple transfer, we would obtain limits with respect to the lifting of the raw data to the big structure $\ma{S}=\IR((t))^{\mathrm{LE}}$ of LE-series. But we want to have a formulation where the limits are taken with respect to the globally subanalytic functions we start with. For this purpose we use the results of the previous Section 4.

\vs{0.5cm}
{\bf 5.2 Theorem} (Lebesgue's theorem on dominated convergence)

\vs{0.1cm}
	{\it Let $f:\IP^{n+1}\to \IP,$ $(s,x)\mapsto f(s,x)=f_s(x),$ be globally subanalytic. Assume that there is some integrable globally subanalytic function $h:\IP^n\to \IP$ such that $|f_s|\leq |h|$ for all sufficiently large $s\in \IP$.
	Then the globally subanalytic function $\lim_{s\to \infty,s\in \IP}f_s$ is integrable and
	\[\int \lim_{s\to \infty,s\in \IP}f_s\,d\lambda_{\IP,n}=\lim_{s\to \infty,s\in \IP}\int f_s\, d\lambda_{\IP,n}.\]}

\vs{0.1cm}
{\bf Proof:}

\vs{0.1cm}
	Note that, by o-minimality and by the assumption that $|f_s|\leq|h|$ for all sufficiently large $s$, $\lim_{s\to \infty,s\in \IP}f_s$ exists and is an integrable globally subanalytic function on $\IP^n$. Note that
	\[\big(\lim_{s\to \infty,s\in\IP} f_s\big)_\ma{S}=\lim_{s\to \infty,s\in\ma{S}} (f_s)_\ma{S}.\]
	Let
	\[F:\IP\to \IP[X], s\mapsto \int f_s\,d\lambda_{\IP,n}.\]
	Then $F$ is a $\IP$-constructible function by Proposition 4.10.
	By the parametric version of Construction 2.6 and by applying the familiar transfer argument, we get that
	\[\int \lim_{s\to \infty,s\in \IP}f_s\,d\lambda_{\IP,n}=\lim_{s\to \infty,s\in\ma{S}}F_{\ma{S}}(s).\]
	Applying the globally subanalytic map $x\to 1/x$ which maps $\infty$ to $0$, we obtain the claim by Theorem 4.20.
\hfill$\Box$

\vs{0.5cm}
{\bf 5.3 Corollary}

\vs{0.1cm}
	{\it Let $A$ be a globally subanalytic subset of $\IP_{\geq 0}\times \IP^{n}$.
	\begin{itemize}
		\item[(1)] (Continuity from below)
		Assume that $A_{s_1}\subset A_{s_2}$ for all $0\leq s_1\leq s_2$ and that there is some globally subanalytic subset $B$ of $\IP^n$ with $\lambda_{\IP,n}(B)<\infty$ such that $A_s\subset B$ for all $s\geq 0$.
		Then
		\[\lim_{s\to \infty,s\in \IP}\lambda_{\IP,n}\big(A_s\big)=\lambda_{\IP,n}\big(\bigcup_{s\geq 0}A_s\big).\]
		\item[(2)] (Continuity from above)
		Assume that $A_{s_1}\supset A_{s_2}$ for all $0\leq s_1\leq s_2$ and that $\lambda_{\IP,n}(A_0)<$ $\infty$.
		Then
		\[\lim_{s\to \infty,s\in \IP}\lambda_{\IP,n}\big(A_s\big)=\lambda_{\IP,n}\big(\bigcap_{s\geq 0}A_s\big).\]
	\end{itemize}}
{\bf Proof:}

\vs{0.1cm}
	Apply Theorem 5.2 and Proposition 2.7(2).
	\hfill$\Box$

\vs{0.5cm}
\noi The continuity from above can be viewed as a substitute for $\sigma$-continuity from above in the usual Lebesgue theory.
The continuity from below can be viewed as a partial substitute for $\sigma$-continuity from below which is equivalent to $\sigma$-additivity.
We need here the additional assumption that the union is contained in a set of finite measure. This assumption is necessary as the following example shows.

\vs{0.5cm}
{\bf 5.4 Example}

\vs{0.1cm}
{\it 
	For $s\geq 1$ let
	\[A_s:=\big\{(x,y)\in \IP^2\mid 1\leq x\leq s, 0\leq y\leq \frac{1}{x}\big\}.\]
	Then $\lambda_{\IP,2}\big(\bigcup_{s\geq 1}A_t\big)=\infty$ but
	$\lambda_{\IP,2}(A_s)=\log s$ does not have a limit in $\IP[X]$ as $s$ tends to $\infty$.}

\vs{0.5cm}
\noi The following example shows that, even in the finite case, the measure is not $\sigma$-additive (cf. the introduction).

\vs{0.5cm}
{\bf 5.5 Example}

\vs{0.1cm}
{\it 
	For $j\in\IN$ let
	$A_j:=[j,t^{-1/j}]$.
	Then $A_j\searrow \emptyset$ but
	$\lambda_{\IP,1}(A_j)$ does not tend to $0$ as $j$ tends to $\infty$ since
	$\lambda_{\IP,1}(A_j)=t^{-1/j}-j$
	for all $j$ by 2.3(5).}

\vs{0.5cm}
{\bf 5.6 Theorem} (Differentiation)

\vs{0.1cm}
	{\it Let $k\in\IN$, let $U$ be an open globally subanalytic subset of $\IP^k$ and let
	$f:U\times \IP^n\to \IP, (s,x)\mapsto f(s,x),$ be globally subanalytic.
	Assume that
	
	\begin{itemize}
		\item[(a)] for all $s\in U$, $f_s:\IP^n\to \IP, x\mapsto f(s,x),$ is integrable,
		\item[(b)] for all $x\in \IP^n$, the partial derivatives of the function $f_x:U\to \IP, s\mapsto f(s,x),$ exist,
		\item[(c)] there is an integrable globally subanalytic function $g:\IP^n\to \IP$ such that 
		\[|\big(\partial f/\partial s_j\big)(s,x)|\leq |g(x)|\] 
		for all  $j\in\{1,\ldots,k\}$ and all $(s,x)\in U\times \IP^n$.
	\end{itemize}
	Then the partial derivatives of the constructible function
	\[\varphi:U\to \IP[X], s\mapsto \int_{\IP^n} f(s,x)\,d\lambda_{\IP,n}(x),\]
	exist and
	\[\frac{\partial \varphi}{\partial s_j}(t)=\int_{\IP^n}\frac{\partial f}{\partial s_j}(s,x)\,d\lambda_{\IP,n}(x)\]
	for all $j\in\{1,\ldots,k\}$ and all $s\in U$.}

\vs{0.1cm}
{\bf Proof:}

\vs{0.1cm}
	Using Theorem 5.2 on dominated convergence, the usual proof (see for example [1, \S 16 I]) can be adjusted to obtain that the partial derivatives of the $\IP$-constructible function $\varphi:U\to \IP[X]$ exist and that the given equalities are valid.
	Note that, by o-minimality, the mean value theorem holds for globally subanalytic functions that are differentiable.
	
\hfill$\Box$

\vs{0.5cm}
{\bf 5.7 Theorem} (Fundamental theorem of calculus)

\vs{0.1cm}
{\it 
	Let $I$ be an open subinterval of $\IP$ and let $f:I\to \IP$ be a globally subanalytic function that is continuous.
	\begin{itemize}
		\item[(1)]
		Let $a\in I$. The function
		\[F:I\to \IP[X], x\mapsto \int_a^xf(s)\,d\lambda_{\IP,1}(s),\]
		is an antiderivative of $f$.
		\item[(2)]
		Let $G$ be an antiderivative of $f$ on $I$. For $a,b\in I$ we have
		\[\int_a^bf(x)\,d\lambda_{\IP,1}=G(b)-G(a).\]
	\end{itemize}}
{\bf Proof:}

\vs{0.1cm}
	(1):
	By Proposition 4.10, $F$ is $\IP$-constructible.
	The usual proof of the fundamental calculus gives that $F$ is differentiable with $F'=f$.
	
	\vs{0.2cm}
	\noi (2): This follows from (1) and Proposition 4.26.
\hfill$\Box$

\vs{0.5cm}
{\bf 5.8 Corollary}

\vs{0.1cm}
	{\it Let  $f:I\to \IP$ be a globally subanalytic function that is continuous.
	Then $f$ has an antiderivative.}

\vs{0.5cm}
{\bf 5.9 Example}

\vs{0.1cm}
	{\it Let $n\in\IZ$. The antiderivative of $x^n$ on $\IP_{>0}$ is given up to additive constants by $x^{n+1}/(n+1)$ if $n\neq -1$ and $\log x$ if $n=-1$.}

\subsection{Fubini's theorem}

\noi So far we cannot establish Fubini's theorem in the above setting. The reason is that one obtains, when integrating with parameters, constructible functions which are not necessarily globally subanalytic. We generalize the construction of Section 2.3, using the results of Cluckers and D. Miller [10, 11, 12], to obtain a version of Fubini's theorem.

\vs{0.5cm}
{\bf 5.10 Construction}	

\vs{0.1cm}	
	Let $f:\IP^n\to \IP$ be $\IP$-constructible.
	We say that $f$ is integrable and define its integral $\int_{\IP^n}f(x)\,dx\in \ma{S}$ as follows.
	Take a formula $\phi(x,s,y)$ in the language $\ma{L}_{\an,\log}$, $x=(x_1,\ldots,x_n), y=(y_1,\ldots,y_q)$, and a point $a\in \IP^q$ such that $\mathrm{graph}(f_\ma{S})=$ $\phi(\ma{S}^{n+1},a)$. We choose thereby $\phi(x,s,y)$ in such a way that $\phi(\IR^{n+1+q})$ is the graph of a constructible function
	$g:\IR^{n+q}\to \IR$.
	By Fact 1.10 there are constructible functions $h:\IR^q\to \IR$ and $F:\IR^q\to \IR$ such that, for every $c\in \IR^q$, the following holds:
	\begin{itemize}
		\item[(1)]
		$g_c:\IR^n\to \IR$ is integrable over $\IR^n$ if and only if $h(c)=0$,
		\item[(2)]
		if $g_c$ is integrable then $\int_{\IR^n}g_c(x)\,dx=F(c)$.
	\end{itemize}
	The graph of the functions $F$ and $h$ are defined in $\IR_{\an,\exp}$ by  $\ma{L}_{\an,\exp}$-formulas $\psi(y,z)$ and $\chi(y,z)$, respectively .
	These formulas define in $\ma{S}$ the graph of a function $F_\ma{S}:\ma{S}^p\to \ma{S}$ and of a function $h_\ma{S}:\ma{S}^p\to\ma{S}$. The values $F_\ma{S}(a)$ and $h_\ma{S}(a)$ do not depend on the choices of $\phi$, $a$ and $\psi,\chi$.
	We say that $f$ is integrable if $h_\ma{S}(a)=0$ and in this case we set
	$\int_{\IP^n}f(x)\,dx:=F_{\ma{S}}(a)$.

\vs{0.5cm}
\noi Note that Theorem 3.11 can be generalized to this setting. We write again $\int_{\IP^n}f\,d\lambda_{\IP,n}$.
Doing the same construction with parameters we obtain the following (compare with Proposition 4.10):

\vs{0.5cm}
{\bf 5.11 Proposition}

\vs{0.1cm}
{\it 	Let $f:\IP^{q+n}\to \IP[X]$ be $\IP$-constructible. The following holds:
	\begin{itemize}
		\item[(1)] There is a $\IP$-constructible function $g:\IP^q\to\IP[X]$ such that
		\[\mathrm{Fin}(f):=\big\{s\in \IP^q\mid f_s\mbox{ is integrable}\big\}\]
		equals the zero set of $g$.
		\item[(2)] There is a $\IP$-constructible function $F:\IP^q\to\IP[X]$ such that
		\[\int_{\IP^n}f_s(x)\,d\lambda_{\IP,n}(x)=F(s)\]
		for all $s\in \mathrm{Fin}(f)$.
	\end{itemize}}

\vs{0.2cm}
\noi Using the usual transfer argument, we obtain Fubini's theorem.

\vs{0.5cm}
{\bf 5.12 Theorem} (Fubini's theorem)

\vs{0.1cm}
{\it 
	Let $f:\IP^{m+n}\to \IP[X]$ be a $\IP$-constructible function that is integrable.
	Let $g:\IP^m\to \IP[X]$ be a $\IP$-constructible function with
	\[g(x)=\int_{\IP^n}f(x,y)\,d\lambda_{\IP,n}(y)\]
	for all $x\in \IP^m$ such that
	$f_x:\IP^n\to\IP[X]$ is integrable.
	Then $g$ is integrable and
	\[\int_{\IP^{m+n}}f(x,y)\,d\lambda_{\IP,m+n}(x,y)=\int_{\IP^m} g(x)\,d\lambda_{\IP,m}(x).\]}

\section{An application}

\noi The Stone-Weierstra\ss $ $ theorem on uniform approximation of continuous functions on bounded and closed intervals by polynomials
does not hold for semialgebraic functions on non-archimedean real closed fields (see [3, Example 8.8.6] for the field of Puiseux series).
By the approximation theorem of Efroymson ([3, Theorem 8.8.4]) it follows that continuous semialgebraic functions can be uniformly approximated by Nash functions. These are functions that are semialgebraic and $C^\infty$.

\noi By integration techniques (smoothing by convolution), we can extend the latter to the
case of $T_\an$-models with archimedean value group. As in the previous sections, we stick to the case of the field of Puiseux series.

\vs{0.2cm}
\noi We consider the semialgebraic function $\Phi:\IP\to\IP, s\mapsto 1/\big(\pi(1+s^2)\big)$.

\vs{0.5cm}
{\bf 6.1 Remark}
	{\it \begin{itemize}
		\item[(1)]
		The function $\Phi$ is integrable and $\int_\IP\Phi(s)\,d\lambda_{\IP,1}(s)=1$.
		\item[(2)]
		Let $r\in \IP_{>0}$.
		Then
		\[\lim_{h\to 0,h\in\IP}\frac{1}{h}\int_{|s|>r}\Phi\big(s/h\big)\,d\lambda_{\IP,1}(s)=0.\]
	\end{itemize}}
{\bf Proof:}

\vs{0.1cm}
	(1): This holds in the real case. We get the claim by Proposition 3.21.
	
	\vs{0.2cm}
	\noi (2): The constructible antiderivative of $\Phi$ is given by the globally subanalytic function $\arctan_\IP$ which is the lifting of the real arctangent to $\IP$.
	Applying the transformation formula Theorem 5.1 and the fundamental theorem of calculus Theorem 5.7 we have that
	\[\frac{1}{h}\int_{|s|>r}\Phi\big(s/h\big)\,d\lambda_{\IP,1}(s)=\int_{|s|>r/h}\Phi(s)\,d\lambda_{\IP,1}(s)=\pi-2\arctan_\IP\big(r/h\big).\]
	Since $\lim_{x\to \infty,x\in \IP}\arctan_\IP(x)=\lim_{x\to \infty,x\in \IR}\arctan(x)=\pi/2$ we obtain the claim.
\hfill$\Box$

\vs{0.5cm}
\noi By Remark 6.1 one can call $\Big(\Phi_h(s)\Big)_{h\in \IP_{>0}}:=\Big(\frac{1}{h}\Phi\big(s/h\big)\Big)_{h\in \IP_{>0}}$ a {\bf Dirac family} on $\IP$.

\vs{0.5cm}
{\bf 6.2 Remark}

\vs{0.1cm}
{\it 
	Let $g:\IP\to \IP$ be a globally subanalytic function that is bounded. Then
	$g(s)\Phi_h(s-x)$ is integrable for all $x\in R$ and all $h\in \IP_{>0}$.}

\vs{0.1cm}
{\bf Proof:}
	
	\vs{0.1cm}
	This follows from Remark 6.1(1).
\hfill$\Box$

\vs{0.5cm}
{\bf 6.3 Definition}

\vs{0.1cm}
	Let $g:\IP\to\IP$ be a globally subanalytic function that is bounded. For $h\in \IP_{>0}$ let
	\[S_hg:\IP\to \IP[X], x\mapsto \int_\IP g(s)\Phi_h(s-x)\,d\lambda_{\IP,1}(s).\]
	The function $S_h(g)$ is the {\bf convolution} of $g$ with $\Phi_h$.

\vs{0.5cm}
\noi We obtain the usual smoothing property:

\vs{0.5cm}
{\bf 6.4 Proposition}

\vs{0.1cm}
{\it 
	Let $g:\IP\to\IP$ be a globally subanalytic function that is continuous. Assume that the support of $g$ is bounded.
	Then the following holds:
	\begin{itemize}
		\item[(1)]
		$S_hg$ is $C^\infty$ for all $h\in \IP_{>0}$.
		\item[(2)]
		For every $\varepsilon\in \IP_{>0}$ there is $h\in \IP_{>0}$ such that $|g(x)-S_hg(x)|<\varepsilon$ for all $x\in \IP$.
	\end{itemize}}
{\bf Proof:}

\vs{0.1cm}
	(1): $\Phi$ is $C^\infty$ and all derivatives of $\Phi$ are bounded on $\IP$. We get the claim by applying Theorem 5.6 repeatedly.
	
	\vs{0.2cm}
	\noi (2): The classical proof (see Bourbaki [5, VIII \S 4]) works in this setting.
\hfill$\Box$

\vs{0.5cm}
{\bf 6.5 Theorem}

\vs{0.1cm}
{\it 
	Let $a,b\in \IP$ with $a<b$ and let $f:[a,b]\to \IP$ be globally subanalytic and continuous.
	Let $\varepsilon\in\IP_{>0}$. Then there is some open interval $I$ in $\IP$ containing $[a,b]$ and some globally subanalytic function $u:I\to \IP$ that is $C^\infty$ such that $|f(x)-u(x)|<\varepsilon$ for all $x\in [a,b]$.}

\vs{0.1cm}
{\bf Proof:}

\vs{0.1cm}
	We may assume that $\varepsilon\in \mathfrak{m}_\IP$. We choose a continuous globally subanalytic function $g:\IP\to\IP$ that extends $f$ and has a bounded support.
	By Proposition 6.4(2) there is some $h\in \IP_{>0}$ such that $|g(x)-S_hg(x)|<\varepsilon$ for all
	$x\in \IP$. Let $v:=S_hg$. Then $v$ is $\IP$-constructible by Proposition 4.10. By Proposition 6.4(1) we have that $v$ is $C^\infty$.
	By Theorem 4.27 there are $N\in\IN$ and globally subanalytic functions $h_0,\ldots,h_N$ that are infinitely often differentiable on an open interval $I$ containing  $[a,b]$ such that $v|_{[a,b]}=\sum_{j=0}^Nh_jX^j$.
	Since $\varepsilon\in\mathfrak{m}_\IP$ we obtain from the transcendence of $X$ that $|h_0(x)-f(x)|<\varepsilon$ for all $x\in[a,b]$.
\hfill$\Box$

\vs{1cm}
\noi \footnotesize{\centerline{\bf References}
	\begin{itemize}
\item[(1)] 
H. Bauer:
Measure and Integration Theory.
De Gruyter, 2001.
\item[(2)]
A. Berarducci and M. Otero:
An additive measure in o-minimal expansions of fields. 
{\it Q. J. Math.} {\bf 55} (2004), no. 4, 411-419. 
\item[(3)]
J. Bochnak, M. Costeand M.-F. Roy:
Real Algebraic Geometry.
Springer, 1998.
\item[(4)]
N. Bourbaki:
Integration. I. Chapters 1-6.
Springer, 2004.
\item[(5)] 
N. Bourbaki:
Integration. II. Chapters 7-9.
Springer, 2004.
\item[(6)]
L. Br\"ocker:
Euler integration and Euler multiplication.
{\it Adv. Geom.} {\bf 5} (2005), no. 1, 145-169.
\item[(7)]
R. Cluckers and M. Edmundo:
Integration of positive constructible functions against Euler characteristic and dimension.
{\it J. Pure Appl. Algebra} {\bf 208} (2007), no. 2, 691-698.
\item[(8)] 
R. Cluckers and F. Loeser:
Constructible motivic functions and motivic integration.
{\it Invent. Math.} {\bf 173} (2008), no. 1, 23-121.
\item[(9)]
R. Cluckers and F. Loeser:
Constructible exponential functions, motivic Fourier transform and transfer principle.
{\it Ann. of Math. (2)} {\bf 171} (2010), no. 2, 1011-1065.
\item[(10)] 
R. Cluckers and D. Miller:
Stability under integration of sums of products of real globally subanalytic functions and their logarithms.
{\it Duke Math. J.} {\bf 156} (2011), no. 2, 311-348.
\item[(11)] 
R. Cluckers and D. Miller:
Loci of integrability, zero loci, and stability under integration for constructible functions on Euclidean space with Lebesgue measure.
{\it Int. Math. Res. Not.} {\bf 2012}, no. 14, 3182-3191.
\item[(12)]
 R. Cluckers and D. Miller:
Lebesgue classes and preparation of real constructible functions.
{\it J. Funct. Anal.} {\bf 264} (2013), no. 7, 1599-1642.
\item[(13)]
R. Cluckers, J. Nicaise and J. Sebag:
Motivic integration and its interactions with model theory and Non-Archimedean geometry, Volume I.
{\it London Math. Soc. Lecture Note Ser.} {\bf 383}, Cambridge Univ. Press, 2011.
\item[(14)]
R. Cluckers, J. Nicaise and J. Sebag:
Motivic integration and its interactions with model theory and Non-Archimedean geometry, Volume II.
{\it London Math. Soc. Lecture Note Ser.} {\bf 384}, Cambridge Univ. Press, 2011.
\item[(15)] 
G. Comte, J.-M. Lion and J.-P. Rolin:
Nature log-analytique du volume des sous-analytiques.
{\it Illinois J. Math.} {\bf 44} (2000), no. 4, 884-888.
\item[(16)] 
O. Costin, P. Ehrlich and H. M. Friedman:
Integration on the surreals: A conjecture of Conway, Kruskal and Norton.
arXiv:1505.02478.
\item[(17)]
 N. Cutland: Nonstandard measure theory and its applications.
{\it Bull. London Math. Soc.} {\bf 15} (1983), no. 6, 529-589.
\item[(18)] 
J. Denef and F. Loeser:
Germs of arcs on singular algebraic varieties and motivic integration.
{\it Invent. Math.} {\bf 135} (1999), no. 1, 201-232.
\item[(19)] 
 L. van den Dries: 
Tame Topology and O-minimal
Structures. {\it London Math. Soc. Lecture Notes Series} {\bf
	248}, Cambridge University Press, 1998.
\item[(20)] 
L. van den Dries: Limit sets in o-minimal structures. In: Proceedings of the RAAG Summer School Lisbon 2003 O-minimal structures.
Cuvillier, 2005, 172-215.
\item[(21)] 
L. van den Dries, A. Macintyre and D. Marker:
The elementary theory of restricted analytic fields with exponentiation.
{\it Annals of Mathematics} {\bf 140} (1994), 183-205.
\item[(22)]
L. van den Dries, A. Macintyre and D. Marker:
Logarithmic-exponential power series.
{\it J. London Math. Soc.} (2) {\bf 56} (1997), no. 3, 417-434.
\item[(23)] 
L. van den Dries, A. Macintyre and D. Marker: Logarithmic-exponential series.
{\it Ann. Pure Appl. Logic} {\bf 111} (2001), no. 1-2, 61-113.
\item[(24)]
L. van den Dries and C. Miller:
Geometric categories and o-minimal structures.
{\it Duke Math. J.} {\bf 84} (1996), no. 2, 497-540.
\item[(25)]
 L. van den Dries and P. Speissegger:
O-minimal preparation theorems.
Model theory and applications, 87-116, Quad. Mat., 11, Aracne, Rome, 2002.
\item[(26)]
A. Fornasiero and E. Vasquez Rifo:
Hausdorff measure on o-minimal structures. 
{\it J. Symbolic Logic} {\bf 77} (2012), no. 2, 631-648. 
\item[(27)]
L. Fuchs:
Teilweise geordnete algebraische Strukturen.
Vandenhoeck \& Ruprecht, 1966.
\item[(28)] 
W. Hodges:
A shorter model theory.
Cambridge University Press, 1997.
\item[(29)] 
E. Hrushovski and D. Kazhdan:
Integration in valued fields.
Algebraic geometry and number theory, 261-405, Progr. Math., {\bf 253}, Birkh\"auser,  2006.
\item[(30)]
E. Hrushovski, Y. Peterzil and A. Pillay:
Groups, measures, and the NIP.
{\it J. Amer. Math. Soc.} {\bf 21} (2008), no. 2, 563-596.
\item[(31)] 
T. Kaiser:
On convergence of integrals in o-minimal structures on archimedean real closed fields.
{\it Ann. Polon. Math.} {\bf 87} (2005), 175-192.
\item[(32)] 
T. Kaiser:
First order tameness of measures.
{\it Ann. Pure Appl. Logic} {\bf 163} (2012), no. 12, 1903-1927.
\item[(33)] 
T. Kaiser:
Integration of semialgebraic functions and integrated Nash functions.
{\it Math. Z.} {\bf 275} (2013), no. 1-2, 349-366.
\item[(34)] 
I. Kaplansky:
Maximal fields with valuations.
{\it Duke Math. J.} {\bf 9} (1942), 303-321.
\item[(35)] 
F.-V. Kuhlmann, S. Kuhlmann and S. Shelah:
Exponentiation in power series fields.
{\it Proc. Amer. Math. Soc.} {\bf 125} (1997), no. 11, 3177-3183.
\item[(36)]
S. Kuhlmann:
Ordered exponential fields.
{\it Fields Institute Monographs}, American Mathematical Society, 2000.
\item[(37)] 
J.-M. Lion and J.-P. Rolin:
Int\'{e}gration des fonctions sous-analytiques et volumes des sous-ensembles sous-analytiques.
{\it Ann. Inst. Fourier (Grenoble)} {\bf 48} (1998), no. 3, 755-767.
\item[(38)]
J. Ma\v{r}\'{i}kov\'{a}:
The structure on the real field generated by the standard part map on an o-minimal expansion of a real closed field.
{\it Israel J. Math.} {\bf 171} (2009), 175-195.
\item[(39)]
J. Ma\v{r}\'{i}kov\'{a} and M. Shiota:
Measuring definable sets in o-minimal fields.
{\it Israel J. Math.} {\bf 209} (2015), 687-714.
\item[(40)] 
S. Prie\ss-Crampe:
Angeordnete Strukturen: Gruppen, K\"orper, projektive Ebenen.
Springer, 1983.
\item[(41)] 
A. Robinson:
Non-standard Analysis. North-Holland Publishing Company, 1966.
\item[(42)]
Y. Yin:
Additive invariants in o-minimal valued fields.
arXiv:1307.0224.
\item[(43)] 
Y. Yomdin and G. Comte:
Tame geometry with application in smooth analysis. Lecture Notes in Mathematics, {\bf 1834}. Springer, 2004.
\end{itemize}}

\vs{0.5cm}
Tobias Kaiser\\
University of Passau\\
Faculty of Computer Science and Mathematics\\
tobias.kaiser@uni-passau.de\\
D-94030 Germany

\end{document}